\documentclass[opre,nonblindrev]{informs3}

\DoubleSpacedXI % Made default 4/4/2014 at request
\usepackage{endnotes}
\let\footnote=\endnote
\let\enotesize=\normalsize
\def\notesname{Endnotes}%
\def\makeenmark{\hbox to1.275em{\theenmark.\enskip\hss}}
\def\enoteformat{\rightskip0pt\leftskip0pt\parindent=1.275em
  \leavevmode\llap{\makeenmark}}
  
\usepackage{algorithm}
\usepackage{algorithmic}
\usepackage{threeparttable}
\usepackage{multirow}
\usepackage{multicol}
\usepackage{setspace}
\usepackage{natbib}
 \bibpunct[, ]{(}{)}{,}{a}{}{,}%
 \def\bibfont{\small}%
\TheoremsNumberedThrough     % Preferred (Theorem 1, Lemma 1, Theorem 2)
\ECRepeatTheorems

\EquationsNumberedThrough    % Default: (1), (2), ...
\MANUSCRIPTNO{} % When the article is logged in and DOI assigned to it,
\begin{document}
%%%%%%%%%%%%%%%%

%\setcounter{page}{1}
%\DOI{10.1287/xxxx.0000.0000}%

\RUNAUTHOR{Ren et al.}

\RUNTITLE{Global Optimization for $K$-Center of One Billion Samples}

% \TITLE{Bundling Information Goods of Decreasing Value}
\TITLE{A Global Optimization Algorithm for $K$-Center Clustering of One Billion Samples}

\ARTICLEAUTHORS{%
\AUTHOR{Jiayang Ren$^1$, Ningning You$^2$, Kaixun Hua$^1$, Chaojie Ji$^3$, Yankai Cao$^1$}
\AFF{$^1$Department of Chemical and Biological Engineering, University of British Columbia, Vancouver, BC, Canada, \EMAIL{rjy12307@mail.ubc.ca}, \EMAIL{kaixun.hua@ubc.ca}, \EMAIL{yankai.cao@ubc.ca}}
\AFF{$^2$Antai College of Economics and Management, Shanghai Jiao Tong University, Shanghai, China, \\
\EMAIL{ningyou@sjtu.edu.cn}}
\AFF{$^3$Department of Mathematics, University of British Columbia, Vancouver, BC, Canada, \\
\EMAIL{chaojiej@math.ubc.ca}}
} % end of the block

\ABSTRACT{%
This paper presents a practical global optimization algorithm for the $K$-center clustering problem, which aims to select $K$ samples as the cluster centers to minimize the maximum within-cluster distance. This algorithm is based on a reduced-space branch and bound scheme and guarantees convergence to the global optimum in a finite number of steps by only branching on the regions of centers. To improve efficiency, we have designed a two-stage decomposable lower bound, the solution of which can be derived in a closed form. In addition, we also propose several acceleration techniques to narrow down the region of centers, including bounds tightening, sample reduction, and parallelization. Extensive studies on synthetic and real-world datasets have demonstrated that our algorithm can solve the $K$-center problems to global optimal within 4 hours for \textbf{ten million samples} in the serial mode and \textbf{one billion samples} in the parallel mode. Moreover, compared with the state-of-the-art heuristic methods, the global optimum obtained by our algorithm can averagely reduce the objective function by 25.8\% on all the synthetic and real-world datasets.

% Enter your abstract
}%

%\KEYWORDS{deterministic inventory theory; infinite linear programming duality;
%  existence of optimal policies; semi-Markov decision process; cyclic schedule}

\HISTORY{This paper is accepted by \textit{Managment Science}. The final published version of this article is available at: https://pubsonline.informs.org/doi/10.1287/mnsc.2023.00218.}

\maketitle
%%%%%%%%%%%%%%%%%%%%%%%%%%%%%%%%%%%%%%%%%%%%%%%%%%%%%%%%%%%%%%%%%%%%%%

% Samples of sectioning (and labeling) in OPRE
% NOTE: (1) \section and \subsection do NOT end with a period
%       (2) \subsubsection and lower need end punctuation
%       (3) capitalization is as shown (title style).
%
%\section{Introduction.}\label{intro} %%1.
%\subsection{Duality and the Classical EOQ Problem.}\label{class-EOQ} %% 1.1.
%\subsection{Outline.}\label{outline1} %% 1.2.
%\subsubsection{Cyclic Schedules for the General Deterministic SMDP.}
%  \label{cyclic-schedules} %% 1.2.1
%\section{Problem Description.}\label{problemdescription} %% 2.

% Text of your paper here

\section{Introduction}
\label{sec: intro}

Cluster analysis is a task to group similar samples into the same cluster while separating less similar samples into different clusters. It is a fundamental unsupervised machine learning task that explores the character of datasets without the need to annotate cluster classes. Clustering plays a vital role in various fields, such as data summarization~\citep{kleindessner2019fair, hesabi2015data}, customer grouping~\citep{aggarwal2004method}, facility location determination~\citep{hansen2009solving}, and etc. 

There are several typical cluster models, including connectivity-based models, centroid-based models, distribution-based models, density-based models, etc. This work focuses on one of the fundamental centroid-based clustering models called the $K$-center problem. The goal of the $K$-center problem is to minimize the maximum within-cluster distance ~\citep{kaufman2009finding}. Specifically, given a dataset with $S$ samples and the desired number of clusters $K$, the $K$-center problem aims to select $K$ samples from the dataset as centers and to minimize the maximum distance from other samples to its closest center. The $K$-center problem is a combinatorial optimization problem that has been widely studied in theoretical computer science \citep{lim_k-center_2005}. Moreover, it has been intensively explored as a symmetric and uncapacitated case of the $p$-center facility location problem in operations research and management science \citep{garcia-diaz_approximation_2019}, where the number of facilities corresponds to the variable $k$ in a standard $K$-center problem.

% Assuming Euclidean distance??? or change to the MINLP form???

Formally, provided a $K$, the objective function of $K$-center problem can be formulated as follows: 
\begin{align} \label{eq:obj}
	\min \limits_{\mu \in X} \max \limits_{s\in \mathcal{S}} \min \limits_{k\in \mathcal{K}} ||x_s-\mu^k||_2^2
\end{align}
where $X=\{x_1,\ldots, x_S \}$ is the dataset with $S$ samples and $A$ attributes, in which $x_s=[x_{s,1}, ..., x_{s,A}]\in \mathbb{R}^{A}$ is the $s$th sample and $x_{s,a}$ is the $a$th attribute of $i$th sample, $s\in\mathcal{S}:=\{1,\cdots,S\}$ is the index set of samples. 
As to the variables related to clusters, $k \in\mathcal{K}:=\{1,\cdots,K\}$ is the index set of clusters, $\mu:=\{\mu^1,\cdots,\mu^K\}$ represents the center set of clusters, $\mu^k=[\mu^k_{1}, ..., \mu^k_{A}]\in \mathbb{R}^{A}$ is the center of $k$th cluster. Here, $\mu$ are the variables to be determined in this problem. We use $\mu \in X$ to denote the ``centers on samples'' constraint in which each cluster's center is restricted to the existing samples.

% where $\mathcal X=\{x_1,\ldots, x_S \} \in\mathbb{R}^{A \times S}$ is the dataset with $S$ samples and $A$ attributes, in which $x_s=[x_{s,1}, ..., x_{s,A}]\in \mathbb{R}^{A}$ is the $s$th sample and $x_{s,a}$ is the $a$th attribute of $i$th sample, $s\in\mathcal{S}:=\{1,\cdots,S\}$ is the index set of the dataset, $k \in\mathcal{K}:=\{1,\cdots,K\}$ is the cluster set, $\mu:=[\mu^1,\cdots,\mu^K]$ represents the center of each cluster. $\mu$ are the variables to be determined. We use $\mu \in X$ to denote the ``centers on samples'' constraint in which .

\subsection{Literature Review}
The $K$-center problem has been shown to be NP-hard \citep{GONZALEZ1985293}, which means that it is unlikely to find an optimal solution in polynomial time unless $P=NP$ ~\citep{garey1979computers}. As a remedy, heuristic algorithms, which aim to find a good but not necessarily optimal solution, are often used to solve the $K$-center problem on large-scale datasets. The study of exact algorithms, which provide an optimal solution but may hardly be terminated in an acceptable time, is restricted to small-scale datasets due to this poor scalability on larger datasets. % Therefore, heuristic algorithms are broadly-studied to obtain a near-optimal solution quickly, while the study of exact algorithms is limited to small-scale datasets (e.g., thousands of samples).

Regarding heuristic algorithms, there are several 2-approximation algorithms that provide a theoretical guarantee of their distance from the optimal solution for the $K$-center problem, but do not provide a guarantee on their running time \citep{plesnik_heuristic_1987, GONZALEZ1985293, dyer_simple_1985, hochbaum_best_1985, cook_combinatorial_1995}. Among these 2-approximation algorithms, Furthest Point First (FPF) algorithm proposed by \cite{GONZALEZ1985293} is known to be the fastest in practice \citep{mihe2005solving}. It works by starting with a randomly selected center and then adding points that are farthest from the existing centers to the center set. Despite their solution quality guarantee, these 2-approximation algorithms may not always provide close-to-optimal solutions in practice \citep{garcia-diaz_approximation_2019}. Another kind of heuristic methods with a polynomial running time but a weaker solution quality guarantee is also intensively studied in the literature \citep{mihe2005solving, garcia-diaz_when_2017}. Besides heuristic methods, there are also metaheuristic methods that do not have a polynomial running time or a solution quality guarantee, but have been shown to provide near-optimal solutions in some cases \citep{mladenovic_solving_2003, pullan_memetic_2008, davidovic_bee_2011}. In sum, none of these algorithms can deterministically guarantee a global optimal solution for the $K$-center problem.

In contrast to the numerous heuristic algorithms, the study of exact algorithms, which provide the optimal solution but no solution time guarantee, is still struggling with small-scale problems (e.g., thousands of samples). Early exact works are inspired by the relationship between $K$-center and \textit{set-covering} problems \citep{minieka_m-center_1970}. \cite{daskin2000new} transferred the $K$-center problem to a \textit{maximal covering} problem, in which the number of covered samples by $K$ centers is maximized. Then, they proposed an iterative binary search scheme to accelerate the solving procedure. \cite{ilhan_efficient_2001} considered iteratively setting a maximum distance and validating if it can cover all the samples. \cite{elloumi2004new} designed a new integer linear programming formulation of the $K$-center problem, then solved this new formulation by leveraging the binary search scheme and linear programming relaxation. These algorithms have been shown to provide practical results on small-scale datasets with up to $1,817$ samples.

Another research direction models the $K$-center problem as a Mixed Integer Programming (MIP) formulation, allowing for the use of the branch and bound technique to find an optimal solution. However, the vanilla implementations of the branch and bound technique are confined to small-scale datasets with fewer than $250$ samples \citep{brusco_branch-and-bound_2005}. Hence, constraint programming is introduced to address the larger scale $K$-center problems. \citet{hutchison_declarative_2013} designed two sets of variables describing the cluster centers and sample belongings, then updated the solution through constraint propagation and branching. They further reduced the sets of variables and proposed a more general framework in \citet{duong2017constrained}. By involving constraint programming, their works can solve the datasets with up to $5,000$ samples.

Recently, researchers have explored iterative techniques to solve the $K$-center problem on large datasets by breaking it down into smaller subproblems, such as iterative sampling \citep{aloise2018sampling} and row generation \citep{contardo2019scalable}. In \citet{aloise2018sampling}, a sampling-based algorithm was proposed that alternates between an exact procedure on a small subset of the data and a heuristic procedure to test the optimality of the current solution. This algorithm is capable to solve a dataset containing $581,012$ samples within $4$ hours. However, a report about the optimality gap is absent, which is an important measure of solution quality. According to that computing the covering set for a subset of all samples is cheaper than all \citep{chen1987relaxation, chen2009new}, the same research group proposed a row generation algorithm that relies on computing a much smaller sub-matrix \citep{contardo2019scalable}. This approach is able to solve a dataset with 1 million samples to a $6\%$ gap in $9$ hours. However, neither of these methods provides a finite-step convergence guarantee, which results in that they may not always converge to an arbitrarily small gap within a finite number of steps. Therefore, these methods can lead to a nontrivial optimality gap, especially for large datasets. % In \citet{aloise2018sampling}, a sampling-based algorithm was proposed for the $K$-center problem. This sampling algorithm alternates between an exact procedure on a small subset to obtain a solution and a heuristic procedure to test the optimality of the current solution. They reported a numerical result on a tremendously larger dataset containing 581,012 samples within 4 hours. However, their work did not report the optimality gap, an important index to evaluate the solution quality. Inspired by the observation that computing the covering set for a subset of all samples is cheaper than all \citep{chen1987relaxation, chen2009new}, the same research group proposed a row generation algorithm that relies on computing a much smaller sub-matrix. This work can solve the dataset with 1 million samples to a 6\% gap in 9 hours \citep{contardo2019scalable}. However, their works do not provide a finite-step convergence guarantee that the algorithm will converge to an arbitrarily small gap with a finite number of steps. Therefore, these methods often lead to a nontrivial optimality gap, especially for large datasets. 

\subsection{Main contributions}
Recently, \cite{cao_scalable_2019} proposed a reduced-space spatial branch and bound (BB) scheme for two-stage stochastic nonlinear programs. \cite{hua_scalable_2021} adopted this reduced-space BB scheme and Lagrangian decomposition to solve the $K$-means clustering problem with a global optimal guarantee. They solve the large-scale $K$-means problems up to $210,000$ samples to $2.6\%$ optimality gap within $4$ hours. However, these works can not be directly applied to the $K$-center problem. The challenge is that the $K$-center problem minimizes the maximum within-cluster distance instead of the average within-cluster distance. Therefore, utilizing the Lagrangian decomposition method to compute the lower bound is impossible. Moreover, because of the ``centers on samples'' constraint in the $K$-center problem, the direct application of Hua's algorithm will lead to infeasible solutions. 

To address these challenges, we propose a tailored reduced-space branch and bound algorithm for the $K$-center problem. We also design several bounds tightening (BT) and sample reduction methods to accelerate the BB procedure. Our algorithm is unique in that it only branches on the region of centers, which allows us to guarantee convergence to the global optimum within a finite number of steps. In contrast, traditional branch and bound algorithms must branch on all integer variables, which can become computationally infeasible for large-scale problems. By focusing on the limited region of centers, our algorithm is capable to solve even large-scale $K$-center problems.

Specifically, the main contributions of this paper are as follows:
\begin{itemize}
\item  We propose an exact global optimization algorithm based on a tailored reduced-space branch and bound scheme for the $K$-center problem. To increase efficiency, we develop a two-stage decomposable lower bounding method with a closed-form solution, eliminating the need for using any MIP solver in the optimization process. Moreover, the convergence of our algorithm to the global optimum is guaranteed by branching only on the region of centers. % We propose an exact global optimization algorithm based on a tailored reduced-space branch and bound scheme for the $K$-center problem. To improve efficiency, we have designed a two-stage decomposable lower bound method with a closed-form solution, which allows us to avoid using MIP solvers to solve optimization sub-problems. Moreover, the convergence of our algorithm to the global optimum is guaranteed by branching only on the region of centers.

\item We demonstrate that the assignment of clusters can be determined for many samples without knowing the optimal solution. Based on this characteristic, we propose several bounds tightening and sample reduction techniques to further reduce the search space and accelerate the solving procedure. Moreover, we also implement a sample-level parallelization strategy to fully utilize computational resources.

\item An open-source Julia implementation of the algorithm is provided. Extensive studies on $5$ synthetic and $33$ real-world datasets have demonstrated that we can obtain the global solution for datasets with up to 1 billion samples and $12$ features, a feat that has not been achieved so far. Especially, compared with the heuristic methods, the global optimum obtained by our algorithm can averagely reduce the objective function by 25.8\% on all the synthetic and real-world datasets.
\end{itemize}

This paper is an expanded version of our proceeding publication \citep{shi_global_2022} that includes one new acceleration technique called sample reduction and a parallel implementation. These improvements have significantly increased the scale of the optimally solvable $K$-center problem from $14$ million samples to $1$ billion. In this version, we  provide more detailed proof of the global optimum convergence of our algorithm. In addition, we have designed more comprehensive numerical experiments on a broader range of datasets and parameters.
% The current paper is an expanded version of our conference publication \citep{shi_global_2022}. Compared with the conference version, we provide one new accelerate technique called sample reduction and parallel algorithm implementation. These two improvements have extended the optimally solvable $K$-center problem from 14 million samples to 1 billion. Moreover, we also provide more comprehensive numerical experiments in a broader range of datasets and parameters. Finally, we present a more detailed proof of the global optimum convergence of our algorithm in this expanded version.

\subsection{Outline}
This paper is organized as follows: Section 2 introduces a two-stage formulation and a Mixed Integer Nonlinear Programming (MINLP) formulation for the $K$-center problem. Section 3 presents the details of the reduced-space branch and bound algorithm, including the lower bound, upper bound methods, and convergence analysis. Section 4 discusses the accelerating techniques for our BB algorithm, including bounds tightening, sample reduction, and parallel implementation techniques. Section 5 presents the detailed proof of convergence to the global optimum in the finite steps. Section 6 gives extensive numerical results compared with other algorithms. Finally, Section 7 concludes the paper.

\section{$K$-center Formulation}\label{sec:formulation}
\subsection{Two-stage Formulation}
To introduce the lower bounding method in the branch and bound scheme, we first propose a two-stage optimization form of the $K$-center Problem \ref{eq:obj}. The first-stage problem is as follows:
\begin{equation}\label{eqn:two}
    z =  \min_{\mu \in X\cap M_0} \max_{s\in \mathcal{S} }Q_s(\mu).
\end{equation}
where the center set $\mu$ is the so-called first-stage variable, $Q_s(\mu)$ is the optimal value of the second-stage optimization problem:
\begin{equation}
    \begin{aligned}
          Q_s(\mu) = \min \limits_{k\in \mathcal{K}}||x_s-\mu^k||_2^2
    \end{aligned}
\end{equation}
We denote a closed set $M_0:=\{\mu\,\mid\,\underbar{$\mu$} \leq \mu \leq \bar{\mu}\}$ as the region of centers, where $\underbar{$\mu$}$ is the lower bound of centers and $\bar{\mu}$ is the upper bound, i.e., $\underbar{$\mu$}^k_a=\min \limits_{s\in \mathcal{S}} X_{s,a}$, $\bar{\mu}^k_a=\max \limits_{s\in \mathcal{S}} X_{s,a}$, $\forall k \in \mathcal{K}$, $a\in\{1,\cdots, A\}$. Here, the constraint $\mu \in M_0$ is introduced to simplify the discussion of the BB scheme. Since $M_0$ can be inferred directly from data, it will not affect the optimal solution of Problem \ref{eq:obj}. Constraint $\mu \in X\cap M_0$ means the center of each cluster is selected from the samples belonging to the intersection set of the corresponding region $M_0$ and the dataset $X$

\subsection{MINLP Formulation}
To introduce the bounds tightening and sample reduction methods, we propose a MINLP formulation of the $K$-center Problem \ref{eq:obj}:
\begin{subequations}\label{eqn:overall}
    \begin{align}
     \min  \limits_{\mu, d, b, \lambda} &\;  d_{*}     \label{eqn:overall:obj}   \\ 
       \rm{s.t.} \;\; & d_s^k \geq ||x_s-\mu^k||_2^2   \label{eqn:overall:dis}    \\
       & -N_1(1-b_s^k) \leq  d_s^*-d_s^k \leq 0 \label{eqn:overall:bigM:b}   \\ 
       & d_*\geq d_s^* \label{eqn:overall:d}\\
        & \sum_{k \in\mathcal{K}} b_s^k = 1 \label{eqn:overall:unique:b} \\
        & b_s^k\in\{0, 1\} \\
        & -N_2(1-\lambda_s^k)\leq x_s-\mu^k \leq N_2(1-\lambda_s^k) \label{eqn:overall:bigM:lambda}\\
        & \sum_{s \in\mathcal{S}} \lambda_s^k = 1 \label{eqn:overall:unique:lambda} \\
        & \lambda_s^k \in \{0,1\}\\
        & b_s^k\geq \lambda_s^k \label{eqn:overall:logic} \\
        & s\in\mathcal{S}, k \in\mathcal{K}
    \end{align}
\end{subequations}
where $d_s^k$ represents the distance between sample $x_s$ and center $\mu^k$, $d_s^*$ denotes the distance between $x_s$ and the center of its cluster, $N_1$ and $N_2$ are both arbitrary large values. $b_s^k$ and $\lambda_s^k$ are two binary variables. $b_s^k$ is equal to 1 if sample $x_s$ belongs to the $K$th cluster, and 0 otherwise. $\lambda_s^k$ is equal to 1 if $x_s$ is the center of the $K$th cluster $\mu^{k}$, and 0 otherwise. 

Constraint \ref{eqn:overall:bigM:b} is a big M formulation and ensures that $d_s^*=d_s^k$ if $b_s^k=1$ and $d_s^*\leq d_s^k$ otherwise. Constraint \ref{eqn:overall:unique:b}  guarantees that sample $x_s$ belongs to one cluster. We also adopt Constraint \ref{eqn:overall:bigM:lambda}, \ref{eqn:overall:unique:lambda} and \ref{eqn:overall:logic} to represent the ``centers on samples'' constraints, $\mu \in X$. Specifically, Constraint \ref{eqn:overall:bigM:lambda} uses a big M formula to make sure that $\mu^k =x_s$ if $\lambda_s^k=1$ and Constraint \ref{eqn:overall:unique:lambda} confirms that each center can only be selected on one sample. Constraint \ref{eqn:overall:logic} ensures that if $x_s$ is the center of the $K$th cluster, then it is assigned to the $K$th cluster. It should be noted that the global optimizer CPLEX also relies on this formulation to solve the $K$-center problem.

\section{Tailored Reduced-space Branch and Bound Scheme}\label{sec:bb}
This section introduces a tailored reduced-space branch and bound algorithm for the $K$-center problem with lower and upper bounding methods.

\subsection{Lower Bounds}\label{sec:clt_lb}

In this section, we adopt the two-stage formulation and derive a closed-form solution to obtain the lower bound of the $K$-center Problem \ref{eq:obj}.

At each node in the BB procedure, we deal with a subset of $M_0$, which is denoted as $M$, and solve the following problem concerning $M$:
\begin{equation}\label{eqn:clt_ndpb}
    z(M) = \min_{ \mu \in X\cap M}\max_{s\in\mathcal{S}} Q_s(\mu)
\end{equation}

This problem can be equivalently reformulated as the following problem by duplicating $\mu$ across samples and enforcing them to be equal:
\begin{subequations}\label{eqn:clt_lift_ndpb}
 \begin{align}
        \min_{\mu_s\in X\cap M} & \max_{s\in\mathcal{S}} Q_s(\mu_s) \\
        \textrm{s.t.} \quad & \mu_s=\mu_{s+1}, s\in \{1,\cdots,S-1\} \label{eqn:non-anticipativity}
  \end{align}
\end{subequations}

We call constraints \ref{eqn:non-anticipativity} the non-anticipativity constraints. By removing the ``centers on samples'' constraint $\mu \in X$ and the non-anticipativity constraints \ref{eqn:non-anticipativity}, we attain a lower bound formulation as follow:
\begin{equation}\label{eqn:lb_pb_minmax}
    \beta(M):= \min_{\mu_s\in M} \max_{s\in\mathcal{S}} Q_s(\mu_s).
\end{equation}

With constraints relaxed, the feasible region of Problem \ref{eqn:lb_pb_minmax} is a superset of Problem \ref{eqn:clt_lift_ndpb}'s feasible region. Therefore, it is obvious that $\beta(M)\leq z(M)$. 

In Problem \ref{eqn:lb_pb_minmax}, since $\mu$ of each sample is independent, it is obvious that:
\begin{equation}\label{eqn:lb_pb}
    \beta(M)=  \max_{s\in\mathcal{S}}\min_{\mu_s\in M} Q_s(\mu_s).
\end{equation}

Clearly, problem \ref{eqn:lb_pb} can be decomposed into $S$ subproblems with $\beta(M)=\max \limits_{s\in \mathcal{S}}\beta_s(M)$:
\begin{equation}\label{eqn:lb_sub}
    \beta_s(M) = \min_{\mu\in M} Q_s(\mu).
\end{equation}

Denote the region of $k$th cluster's center as $M^k:=\{\mu^k: \underbar{$\mu$}^k \leq \mu^k\leq \bar{\mu}^k\}$ where $\underbar{$\mu$}^k$ and $\bar{\mu}^k$ are the lower and upper bound of $\mu^k$ respectively. 
Since $Q_s(\mu) =\min \limits_{k\in \mathcal{K}}||x_s-\mu^k||_2^2$, we have 
\begin{equation}\label{eqn:lb_sub2}
    \beta_s(M) =   \min \limits_{k\in \mathcal{K}}  \min_{\mu^{k}\in M^k}  ||x_s-\mu^{k}||_2^2,
\end{equation}
which can be further decomposed into $K$ subsubproblems with $\beta_s(M) {=} \min\limits_{k\in \mathcal{K}} \beta_{s}^k(M^k)$:
\begin{equation}\label{eqn:lb_subsub}
    \beta_{s}^k(M^k) =\min_{\mu^{k}\in M^k}  ||x_s-\mu^{k}||_2^2.
\end{equation}

The analytical solution to Problem \ref{eqn:lb_subsub} is: ${\mu_{a}^k}^* = \text{mid}\{\underbar{$\mu$}^k_a, \ x_{s,a}, \ \bar{\mu}^k_a\}, \forall a\in \{1,\cdots,A\}$. Consequently, the closed-form solution to Problem \ref{eqn:lb_pb_minmax} can be easily computed by the max-min operation on all the samples.

\subsection{Upper Bounds} \label{sec:clt_ub}
At each node in the BB procedure, the upper bounds of Problem \ref{eqn:clt_ndpb} can be obtained by fixing the centers at a candidate feasible solution $\hat{\mu}\in X\cap M$. In this way, we can compute the upper bound base on the following equation:
\begin{equation}\label{eqn:ub_gp_ct}
    \alpha(M)=\max \limits_{s\in \mathcal{S}}\min \limits_{k\in \mathcal{K}}||x_s-\hat{\mu}^k||_2^2
\end{equation}

Since $\hat{\mu}$ is a feasible solution, we have $z(M)\leq \alpha(M)$, $\forall \hat{\mu}\in X\cap M$. In our implementation, we use two methods to obtain the candidate feasible solutions. At the root node, we use a heuristic method called Farthest First Traversal \citep{GONZALEZ1985293} to obtain a candidate solution $\hat{\mu}\in X\cap M_0$. Using this method, we randomly pick an initial point and select each following point as far as possible from the previously selected points. Algorithm \ref{alg: fft} describes the details of the farthest first traversal, where $d(x_s, T)$ represents the minimum distance from sample $x_s$ to any sample in set $T$. We use $FFT(M_0)$ to denote the upper bound obtained using this approach. At a child node with center region $M$, for each cluster, we select the data sample closest to the middle point of $M^k$ as $\hat{\mu}^k$, and obtain the corresponding upper bound $\alpha(M)$. 

\subsection{Branching}

Our algorithm only needs to branch on the region of centers, $M:=\{\mu: \underbar{$\mu$} \leq \mu\leq \bar{\mu}\}$, to guarantee convergence, which would be theoretically discussed in Section 5, o. Since the desired number of clusters is $K$ and the number of attributes is $A$, the number of possible branching variables is $K\times A$. The selection of branching variables and values will dramatically influence the BB procedure's efficiency. In our implementation, we select the max-range variable at each node as the branching variable and the midpoint of this variable as the branching value. %Specifically, the branching variable is selected by $\argmax\limits_{k,a} |\bar{\mu}^k_a-\underbar{$\mu$}^k_a|$. Suppose the index of the branching variable is $k,a$, the branching value is $\frac{\bar{\mu}^k_a-\underbar{$\mu$}^k_a}{2}$. This max-range branching method can help to balance the search space and reduce computational costs.

\subsection{Branch and Bound Scheme}\label{sec:bb_clst}
The detailed reduced-space branch and bound algorithm for the $K$-center Problem \ref{eq:obj} are given in the Algorithm \ref{alg: bb_sche}. In the algorithm, We use $relint(.)$ to denote the relative interior of a set.  We can also establish the convergence of the branch-and-bound scheme in Algorithm \ref{alg: bb_sche}. The BB procedure can generate a monotonically non-ascending sequence $\{\alpha_i\}$ and a monotonically non-descending sequence $\{\beta_i\}$. We can show that they both converge to $z$ in a finite number of steps.

\begin{theorem}\label{theorem: conv_finite} 
\textit{Algorithm \ref{alg: bb_sche} is convergent to the global optimal solution after a finite step $L$, with $\beta_L=z=\alpha_L$, by only branching on the region of centers.}
\end{theorem}

Since the following acceleration techniques also influence the global convergence in Section \ref{sec:acc}. We present the detailed proof of Theorem \ref{theorem: conv_finite} in Section \ref{sec:convergence} after introducing the acceleration techniques.

\begin{minipage}[!htp]{\textwidth}
\begin{multicols}{2}%[3in]{\dimexpr.5\textwidth}
 \begin{minipage}[T]{0.49 \textwidth}
\begin{algorithm}[H]
   \caption{Branch and Bound Scheme}
   \scriptsize
   \label{alg: bb_sche}
   \renewcommand{\arraystretch}{0.1}
    \begin{algorithmic}
        \STATE {\bfseries Initialization} 
        \STATE Initialize the iteration index $i\leftarrow 0$;
        \STATE Set $\mathbb{M}\leftarrow\{M_0\}$, and tolerance $\epsilon > 0$;
        \STATE Compute initial lower and upper bounds $\beta_i = \beta(M_0)$, $\alpha_i = FFT(M_0)$ // Alg. \ref{alg: fft} ;
        \STATE Select $K$ farthest initial seeds // Sec.\ref{sec:fbbt_fft};
        \WHILE{$\mathbb{M}\neq \emptyset$}
            \STATE {\bfseries Node Selection}
            \STATE Select a set $M$ satisfying $\beta(M)=\beta_i$ from $\mathbb{M}$ and delete it from $\mathbb{M}$;
            \STATE Update $i\leftarrow i+1$;
            
            \STATE {\bfseries Bounds Tightening}
            \STATE Cluster Assignment // Alg. \ref{alg: assignment};
            \STATE Bounds Tightening // Alg. \ref{alg: fbbt};
            \STATE Obtain the tightened node $\hat{M}$;
            %\STATE Sample Reduction // Alg. \ref{alg: reduction};
            \STATE If $i\ \% \ i_{sr}=0$, Sample Reduction // Alg. \ref{alg: reduction};
            \IF{$\exists |X\cap M^k|>1, k\in \mathcal{K}$}
            \STATE {\bfseries Branching}
            \STATE Find two subsets $M_1$ and $M_2$ s.t. $relint(M_1)\cap relint(M_2) = \emptyset$ and $M_1\cup M_2=M$;
            \STATE Update $\mathbb{M}\leftarrow \mathbb{M}\cup \{M_i\}$, if $ \ X\cap M_i^k \neq \emptyset, \forall k\in\mathcal{K}, i\in{1,2}$;
            \ENDIF
          \STATE {\bfseries Bounding}
            \STATE Compute upper and lower bound $\alpha(M_1)$, $\beta(M_1)$, $\alpha(M_2)$, $\beta(M_2)$;
            \STATE Let $\beta_i\leftarrow \min\{\beta(M')\,\mid\,M'\in\mathbb{M}\}$;
            \STATE Let $\alpha_i\leftarrow \min\{\alpha_{i-1}, \alpha(M_1), \alpha(M_2)\}$;
            \STATE Remove all $M'$ from $\mathbb{M}$ if $\beta(M')\geq\alpha_i$;
            \STATE If $\beta_i-\alpha_i\leq\epsilon$, STOP;
     \ENDWHILE   
    \end{algorithmic}
\end{algorithm}

\vskip -0.4in
\begin{algorithm}[H]
   \caption{Farthest First Traversal}
   \scriptsize
   \label{alg: fft}
    \begin{algorithmic}
        \STATE {\bfseries Initialization} 
        \STATE Randomly pick $s\in \mathcal{S}$;
        \STATE Denote $T$ as the set of $K$ points selected by farthest first traversal;
        \STATE Set $T\leftarrow \{x_s\}$;
        \WHILE{$|T|<K$}
            \STATE Compute $x_s \in \arg \max\limits_{x_s\in X}d(x_s, T)$ to find $x_s$ which is the farthest away from set $T$;
            \STATE $T\leftarrow T\cup \{x_s\}$;
         \ENDWHILE   
    \end{algorithmic}
\end{algorithm}

\end{minipage}

\begin{minipage}[T]{0.49 \textwidth}

\begin{algorithm}[H]
   \scriptsize
   \caption{Cluster Assignment}
   \label{alg: assignment}
\begin{algorithmic}
    \STATE {\textbf{Center Based Assignment}}
    \FOR{sample $x_s\in X$}
        \IF{$b_s^{k}==0, \forall k \in \mathcal{K}$}
            \IF{$\beta_{s}^{k}(M^{k})>\alpha, \forall k \in \mathcal{K} \setminus \{k'\}$}
                \STATE $x_s$ is assigned to cluster $k'$ with $b_s^{k'}=1$;
            \ENDIF
        \ENDIF
    \ENDFOR
    
    \STATE {\textbf{Sample Based Assignment}}
    \IF{All clusters have at least one sample assigned}
        \FOR{sample $x_s\in X$}
            \IF{$\forall k \in \mathcal{K} \setminus \{k'\}$, $\exists\ x_j$ assigned to kth cluster, $||x_s - x_j||^2_2>4\alpha$}
                \STATE $x_s$ is assigned to cluster $k'$ with $b_s^{k'}=1$.
            \ENDIF
        \ENDFOR
    \ENDIF
\end{algorithmic}
\end{algorithm}

\vskip -0.36in
\begin{algorithm}[H]
   \scriptsize
   \caption{Bounds Tightening}
   \label{alg: fbbt}
\begin{algorithmic}
    \STATE Given the current center region $M$ and upper bound $\alpha$
    \FOR{Cluster $k \in \mathcal{K}$}
    \STATE Obtain the assigned sample set $\mathcal{J}^k$ using Alg.\ref{alg: assignment};
    \STATE Compute the ball-based or box-boxed area of each assigned sample, $B_{\alpha}(x_j)$ or $R_{\alpha}(x_j)$;
    \STATE Tighten the center region by $M^k\cap B_{\alpha}(x_j)$ or $M^k\cap R_{\alpha}(x_j)$ , $\forall j \in \mathcal{J}^k$;
    \STATE Further tighten according to the ``centers on samples'' constraint;
    \ENDFOR
 \end{algorithmic}
\end{algorithm}

\vskip -0.36in
\begin{algorithm}[H]
   \scriptsize
   \caption{Sample Reduction}
   \label{alg: reduction}
\begin{algorithmic}
\STATE Initialize the index set of redundant samples as $\mathcal{R}\gets \mathcal{S}$
\FOR{all BB nodes} 
\STATE Obtain the index set of redundant samples for lower bounds, $\mathcal{R}_{LB}$, according to the criterion in Sec. \ref{sec: re_lb};
\STATE Obtain the index set of redundant samples for upper bounds, $\mathcal{R}_{UB}$, according to the criterion in Sec. \ref{sec: re_ub};
\STATE Update the redundant index set, $\mathcal{R}\gets \mathcal{R}\cap \mathcal{R}_{LB}\cap \mathcal{R}_{UB}$;
\ENDFOR
\STATE Delete samples in the redundant set $\mathcal{R}$ from the current dataset.
\end{algorithmic}
\end{algorithm}

\end{minipage}
\end{multicols}
\end{minipage}
\vskip 0.4in

\section{Acceleration Techniques} \label{sec:acc}
Although the lower bound introduced in Section \ref{sec:clt_lb} is enough to guarantee convergence, it might not be very tight, leading to tremendous iterations. Therefore, we propose several acceleration techniques to reduce the search space and speed up the BB procedure. Since Algorithm \ref{alg: bb_sche} only branches on the region of centers $M:=\{\mu: \underbar{$\mu$} \leq \mu\leq \bar{\mu}\}$, we focus on reducing the region of centers to accelerate the solution process while not excluding the optimal solution of the original $K$-center problem.

\subsection{Bounds Tightening Techniques} \label{sec:BT}
In each node, the assignment of many samples (i.e., which cluster the sample is assigned to) can be pre-determined by the geometrical relationship of samples and regions of centers. This information can be further used to reduce the region of $\mu$. 

\subsubsection{Cluster Assignment}\label{sec:fbbt_assign}
The task of cluster assignment is to pre-determine some values of $b_s^k$ in the MINLP Formulation \ref{eqn:overall} at each BB node before finding the global optimal solution.

We first demonstrate the relations between samples and centers. Denote $\alpha$ as the upper bound obtained using methods described in Section \ref{sec:clt_ub}. Then based on Objective \ref{eqn:overall:obj} and Constraint \ref{eqn:overall:d}, we have $d_s^* \leq d_* \leq \alpha$. From Constraint \ref{eqn:overall:dis} and \ref{eqn:overall:bigM:b}, we can conclude that if $b_s^k=1$, then $||x_s-\mu^k||_2^2\leq d_s^* \leq \alpha$. Therefore, we can derive Lemma \ref{lemma:fbbt_1a}:
\begin{lemma}\label{lemma:fbbt_1a}
     If sample $x_s$ is in the $k$th cluster, then $||x_s-\mu^k||_2^2\leq \alpha$, where $\alpha$ is an upper bound of the $K$-center problem.
\end{lemma}

Besides the relation between samples and centers, cluster assignments may also be determined from the distance of two samples. Suppose sample $x_i$ and $x_j$ belong to the $k$th cluster, then from Lemma \ref{lemma:fbbt_1a} we have $||x_i-\mu^k||_2^2\leq \alpha$ and $||x_j-\mu^k||_2^2\leq \alpha$. Thus $||x_i-x_j||_2^2 = ||x_i-\mu^k + \mu^k - x_j||_2^2 \leq  (||x_i-\mu^k||_2 + ||\mu^k - x_j||_2)^2 \leq 4\alpha$. Therefore, we have Lemma \ref{lemma:fbbt_seeds}:

\begin{lemma}\label{lemma:fbbt_seeds}
    If two samples $x_i$ and $x_j$ are in the same cluster, then $||x_i-x_j||_2^2 \leq 4\alpha$ where $\alpha$ is an upper bound of the $K$-center problem.
\end{lemma}

\begin{figure}[tbh]
\begin{center}
\centerline{\includegraphics[width=3.6in]{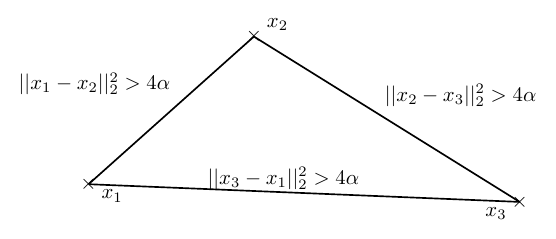}}
\caption{Initial seeds with 3 clusters. In this example, $||x_1-x_2||^2_2>4\alpha$, $||x_2-x_3||^2_2>4\alpha$ and $||x_3-x_1||^2_2>4\alpha$. Therefore, we can arbitrarily assign $x_1, x_2, x_3$ to  3 distinct clusters.}
\label{fig:k_farthest_points}
\vskip -0.2in
\end{center}
\end{figure}

We propose three methods for pre-assigning samples based on these two Lemmas:

\textbf{$K$ Farthest Initial Seeds: }\label{sec:fbbt_fft}
From Lemma \ref{lemma:fbbt_seeds}, if $||x_i-x_j||_2^2 > 4\alpha$, then $x_i$ and $x_j$ are not in the same cluster. At the root node, if we can find $K$ samples with the distance between any two of these samples $x_i$ and $x_j$ satisfying $||x_i-x_j||_2^2 > 4\alpha$, then we can conclude that these $K$ samples must belong to $K$ distinct clusters. Figure \ref{fig:k_farthest_points} shows an example of this property, in which three samples are pre-assigned to 3 distinct clusters. We call these $K$ points initial seeds. To find the initial seeds, every two samples must be as far as possible. Therefore, in our implementation, we use the heuristic Farthest First Traversal (FFT) (Algorithm \ref{alg: fft}) to obtain $K$ farthest points. For about half of the case studies shown in Section \ref{sec:results}, we can obtain the initial seeds using FFT. However, for other cases, initial seeds can not be obtained using FFT, or the initial seeds may not even exist. 

\begin{figure}[tbh]
\begin{center}
\centerline{\includegraphics[width=3.6in]{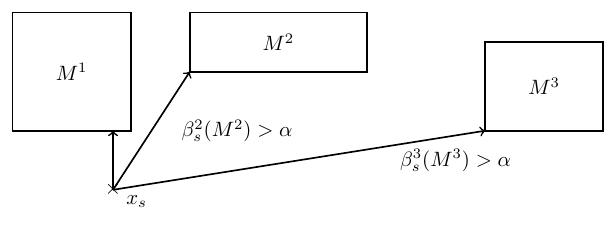}}
\caption{Center-based assignment with 3 clusters. In this example, $\beta_{s}^2(M^2)>\alpha$ ($b_s^2=0$) and $\beta_{s}^3(M^3)>\alpha$ ($b_s^3=0$). Therefore, we assign $x_s$ to the first cluster ($b_s^1=1$).}
\label{fig:center_based_assign}
\vskip -0.2in
\end{center}
\end{figure}

\textbf{Center-Based Assignment:}\label{sec:fbbt_center_based}
From Lemma \ref{lemma:fbbt_1a}, if $||x_s-\mu^k||_2^2>\alpha$, then $x_s$ does not belong to $k$th cluster, which is $b_s^k=0$. Consequently, if we can determine that $b_s^k=0, \forall k \in \mathcal{K} \setminus \{k'\}$, then $b_s^{k'}=1$. However, the value of $\mu$ here is unknown before obtaining the optimal solution. One observation is that if the BB node with region $M$ contains the optimal solution, then we have $\beta_{s}^k(M^k) =\min\limits_{\mu^{k}\in M^k} ||x_s-\mu^{k}||_2^2 \leq ||x_s-\mu^{k}||_2^2$. Therefore, if $\beta_{s}^k(M^k)>\alpha$, sample $x_s$ is not in the $k$th cluster and $b_s^k=0$. In summary, for sample $x_s$, if $\forall k \in \mathcal{K} \setminus \{k'\}$, $\beta_{s}^{k}(M^{k})>\alpha$, then $x_s$ is assigned to cluster $k'$ with $b_s^{k'}=1$. Figure \ref{fig:center_based_assign} illustrates an example in two-dimensional space with three clusters.

This center-based method can be adopted at every node of the BB scheme. Since $\beta_{s}^k(M^k)$ is already obtained when computing the lower bound in Section \ref{sec: re_lb}, there is no additional computational cost. Nevertheless, we do not need to apply this method at the root node since $M_0^1=\cdots=M_0^K$. As the BB scheme continues branching on the regions of centers, $M^k$ becomes more and more different from others. Then more samples can be pre-assigned using this center-based method.

\begin{figure}[tbh]
\vskip -0.2in
\begin{center}
\centerline{\includegraphics[width=3.6in]{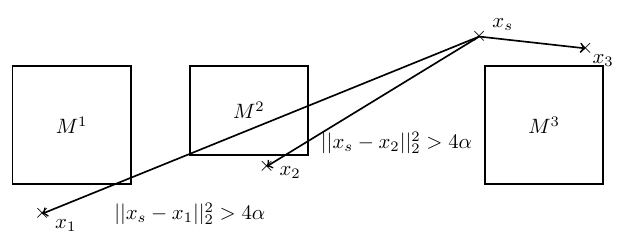}}
\caption{Sample-based assignment with 3 clusters. Assume we already know that $x_1, x_2, x_3$ belong to cluster $1,2$ and $3$, respectively. $x_s$ is the sample to be determined. In this example, $||x_s-x_1||^2_2>4\alpha$ ($b_s^1=0$) and $||x_s-x_2||^2_2>4\alpha$ ($b_s^2=0$). Therefore, $x_s$ is assigned to cluster 3 ($b_s^3=1$).}
\label{fig:sample_based_assign}
\end{center}
\vskip -0.2in
\end{figure}

\textbf{Sample-Based Assignment:}\label{sec:fbbt_sample_based}
Besides utilizing centers to pre-assign samples, assigned samples can also help pre-assign other samples. From Lemma \ref{lemma:fbbt_seeds}, if $||x_i-x_j||_2^2 > 4\alpha$, then $x_i$ and $x_j$ are not in the same cluster. If $x_j$ belongs to $k$th cluster, then obviously $x_i$ cannot be assigned to $k$the cluster and $b_i^k=0$. With this relationship, if all the other $K-1$ clusters are excluded, $x_i$ will be assigned to the remaining cluster. Figure \ref{fig:sample_based_assign} shows an example of the sample-based assignment. 

There is a prerequisite to using this sample-based method. For each cluster, there must be at least one sample already assigned to the cluster. Based on this prerequisite, sample-based assignment is utilized only after at least one sample is pre-assigned for each cluster.

\subsubsection{Bounds Tightening} \label{sec:fbbt}
In this subsection, we adopt the Bounds Tightening (BT) technique and the cluster assignment information to reduce the region of $\mu$. 

\begin{figure}[tbh]
\vskip -0.1in
\begin{center}
\centerline{\includegraphics[width=2.4in]{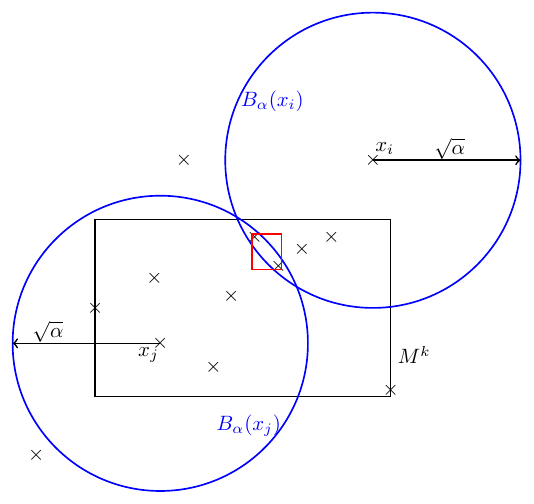}}
\caption{Ball-based bounds tightening in two-dimensional space. In this example, suppose it is determined that two points $x_i$ and $x_j$ belong to the $K$th cluster. We first compute the index set of samples within all balls and original box, $\mathcal{S}^k_{+}(M):= \{s\in \mathcal{S} \ | x_{s}\in X\cap M^k \cap B_{\alpha}(x_i)\cap B_{\alpha}(x_j)\}$. We then generate the smallest box containing these samples in $\mathcal{S}^k_{+}(M)$. The red rectangle is the tightened bounds we obtain.}
\label{fig:ball}
\end{center}
\vskip -0.2in
\end{figure}

\textbf{Ball-based Bounds Tightening: }
For a sample $j$, $B_{\alpha}(x_j){=}\{x  | \ ||x-x_j||_2^2\leq \alpha\}$ represents the ball with center $x_j$ and radius $\sqrt{\alpha}$. By using cluster assignment methods in Section \ref{sec:fbbt_assign}, assuming that sample $j$ belongs to $k$th cluster is already known, by Lemma \ref{lemma:fbbt_1a}, then $\mu^k\in B_{\alpha}(x_j)$ holds. We use $\mathcal{J}^k$ to denote the index of all samples assigned to $k$th cluster, i.e., $\mathcal{J}^k= \{j\in \mathcal{S}\ | \ b_j^k = 1\}$, then $\mu^k\in B_{\alpha}(x_j), \forall j \in \mathcal{J}^k$. Besides this, we also know that $\mu^k \in X\cap M^k$. Denote $\mathcal{S}_{+}^k$ as the index set of samples satisfying all these constraints, $\mathcal{S}^k_{+}(M):= \{s\in \mathcal{S} \ | x_{s}\in X\cap M^k, x_{s} \in B_{\alpha}(x_j), \forall j \in \mathcal{J}^k\}$. In this way, we can obtain a tightened box containing all feasible solutions of $k$th center, $\hat{M}^k {=} \{\mu^k|\hat{\underbar{$\mu$}}^k\leq \mu^k\leq \hat{\bar{\mu}}^k\}$, with the bounds of $a$th attribute in $k$th center to be $\hat{\underbar{$\mu$}}_{a}^k {=} \min\limits_{s\in \mathcal{S}^k_{+}(M)} x^k_{s,a}$ and $\hat{\bar{\mu}}_{s}^k {=} \max\limits_{s\in \mathcal{S}^k_{+}(M)} x^k_{s,a}$. Figure \ref{fig:ball} gives an example of bounds tightening using this method. One challenge of this ball-based bounds tightening method is that it needs to compute the distance of $x_s$ and $x_j$ for all $s\in \mathcal{S}$ and $j \in \mathcal{J}^k$. If we know the assignments of the majority of the samples, we need to do at most $S^2$ times of distance calculation. Note that we only need to do $S*K$ times of distance calculation to compute a lower bound. To reduce the computational time, we set a threshold on the maximum number of balls (default: 50) utilized to tighten bounds in our implementation.  

\begin{figure}[tbh]
\begin{center}
\centerline{\includegraphics[width=2.4in]{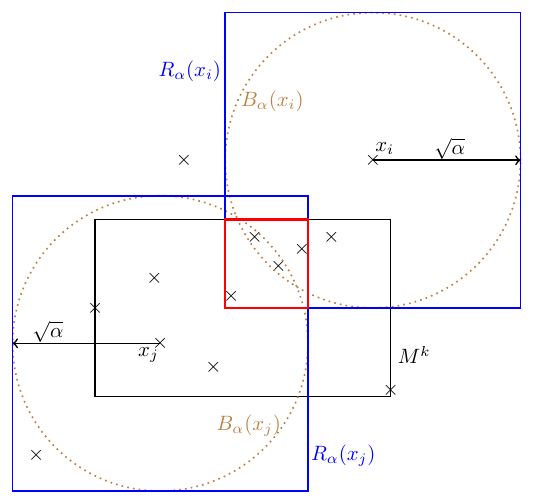}}
\caption{Box-based bounds tightening in two-dimensional space. In this example, we first generate two boxes with $R_{\alpha}(x_i):=\{x| \ x_i-\sqrt{\alpha}\leq x \leq x_i+\sqrt{\alpha}\}$ and $R_{\alpha}(x_j)=\{x|\ x_j-\sqrt{\alpha}\leq x \leq x_j+\sqrt{\alpha} \}$. We then create a tighten bounds with $\hat{M}^k {=} R_{\alpha}(x_i) \cap R_{\alpha}(x_j) \cap M^k$. The red rectangle is the tightened bounds we obtain.}
\label{fig:box}
\end{center}
\vskip -0.3in
\end{figure}

\textbf{Box-based Bounds Tightening:}
Another strategy to reduce the computation burden is based on the relaxation of $B_{\alpha}(x_j)$. For any ball $B_{\alpha}(x_j)$, the closed set $R_{\alpha}(x_j)=\{x \mid x_j-\sqrt{\alpha}\leq x \leq x_j+\sqrt{\alpha} \}$ is the smallest box containing $B_{\alpha}(x_j)$. Then we have $\mu^k\in R_{\alpha}(x_j), \forall j \in \mathcal{J}^k$. Since $R_{\alpha}(x_j)$ and $M^k$ are all boxes, we can easily compute the tighten bounds $\hat{M}^k {=} \bigcap_{j \in \mathcal{J}^k}  R_{\alpha}(x_j) \cap M^k$. Figure \ref{fig:box} gives an example of box-based bounds tightening using this method. Obviously, the bounds generated in Figure \ref{fig:ball} is much tighter, while the method in Figure \ref{fig:box} is much faster. Consequently, if $|\mathcal{J}^k|$ is small for all clusters, the ball-based bounds tightening method gives more satisfactory results. While if $|\mathcal{J}^k|$ is large for any $k$, box-based bounds tightening provides a cheaper alternative. 

\subsubsection{Symmetry Breaking}\label{sec:sb}
Another way to get tighter bounds is based on symmetry-breaking constraints. We add the constraints $\mu_1^1\leq \mu_1^2\leq \cdots \leq \mu_1^K$ in the BB algorithm \ref{alg: bb_sche}, in which $\mu_a^k$ denotes $a$th attribute of $k$th center. Note that symmetry-breaking constraints and FFT-based initial seeds in Section \ref{sec:fbbt_fft} both break symmetry by providing a certain order for the clusters, so they cannot be combined. Our implementation uses symmetric breaking only when initial seeds are not found from FFT at the root node. It should be noted that we also add this symmetry-breaking constraints when using CPLEX to solve the MINLP formulation \ref{eqn:overall} of the $K$-center problem. 

\subsection{Sample Reduction}
Some samples may become redundant during the lower and upper bounding procedure without contributing to the bound improvements. If these samples are proven to be redundant in all the current and future branch nodes, we can conclude they will not influence the bounding results anymore, resulting in sample reduction.

%Hence, only the sample with $\beta_s(M) \geq \beta$ will influence the update of lower bounds. 

\subsubsection{Redundant samples in lower bounding} \label{sec: re_lb}
Denote $\beta$ as the current best lower bound obtained using methods described in Section \ref{sec:clt_lb}. According to Equation \ref{eqn:lb_pb}, lower bound $\beta(M)$ is the maximum value of each sample's optimal value, $\beta_s(M)$. Based on this observation, we further define the best maximum distance of sample $s$ to the center region of $\mu$ as 
\begin{equation}
    \label{eqn: lb_re}
    \alpha_s(M) = \min_{k\in\mathcal{K}}\max_{\mu^k\in M^k}||x_s-\mu^k||^2_2,
\end{equation}
It is obvious that $\beta_s(M) \leq \alpha_s(M)$. If $\alpha_s(M ) < \beta$, we have $\beta_s(M) < \beta$, which means sample $s$ is not the sample corresponding to maximum within-cluster distance. Hence, we can conclude that sample $s$ is a redundant sample in lower bounding for this BB node. 
Moreover, $\forall M' \subset M$, we have $\beta_s(M') \leq \alpha_s(M') \leq \alpha_s(M)$. According to the shrinking nature of center region $M$ and the non-descending nature of lower bound $\beta$, if $\alpha_s(M)< \beta$ is true in a BB node, sample $s$  will remain redundant in all the child nodes of this branch node. It should be noted that $\alpha_s(M)$ can be calculated using an analytical solution similar to $\beta_s(M)$, which is $\mu^k_a= \underbar{$\mu$}^k_a$ if $|\underbar{$\mu$}^k_a - x_{s,a}| > |\bar{\mu}^k_a - x_{s,a}|$, otherwise $\bar{\mu}^k_a$.

\subsubsection{Redundant samples in upper bounding} \label{sec: re_ub}

Obviously, a sample $x_j$ cannot be the center for $k$th cluster if it does not belong to $M^k$. Moreover, according to Lemma \ref{lemma:fbbt_1a}, if a sample $x_j$ is the center for cluster $K$, $||x_i-x_j||^2_2\leq \alpha$ must hold for all the samples $x_i$ assigned to this cluster. Hence, a sample $x_j$ also cannot be the center for $k$th cluster, if there exists a sample $x_i$ assigned to $k$th cluster satisfying $||x_i-x_j||^2_2> \alpha$. If sample $x_j$ cannot be centers for any cluster, we denote this sample $x_j$ as a redundant sample for upper bounding. Since the non-ascending nature of upper bound $\alpha$, if sample $s$ is redundant for upper bounding in a branch node, it will remain redundant in all the child nodes of this branch node. It should be noted that the calculations in this method are identical to Sample-Based Assignment in Section \ref{sec:fbbt_center_based} with no extra calculations introduced in this method.

\subsubsection{Sample reduction} 
If a sample $s$ is redundant in lower bounding, it implies that sample $s$ is not the ``worst-case sample'' corresponding to the maximum within-cluster distance. If a sample $s$ is redundant in upper bounding, then it means that sample $s$ cannot be a center for any cluster. If the sample $s$ is redundant in both lower bounding and upper bounding, then removing this sample will not affect the solution of this BB node and all its child BB nodes. Algorithm \ref{alg: reduction} describes the procedure of sample reduction: first, obtain the redundant samples for lower and upper bounding in each branch node; then, we can delete the samples that are redundant for both lower and upper bounding in all the branch nodes. In our implementation, this sample reduction method is executed for every $i_{sr}$ iterations.  

\subsubsection{Effects on computation}
Sample reduction can reduce the number of samples that need to be explored by deleting redundant samples every $i_{sr}$ iterations, as described in Algorithm \ref{alg: reduction}. It can also accelerate the calculation of lower bounds and bounds tightening at each iteration. For the lower bounding method in Section \ref{sec:clt_lb}, we only need to solve the second-stage problems for  non-redundant samples that have been validated by the lower-bounding criterion in Section \ref{sec: re_lb}. Additionally, once a sample is deemed redundant for lower bounding in a particular node, it will remain redundant in all child nodes of that node. This means that we do not need to solve the second-stage problem for this sample in the current node or any of its child nodes. For the bounds tightening methods in Section \ref{sec:fbbt}, we only need to calculate the bounds based on non-redundant samples that have been validated by the upper-bounding criterion in Section \ref{sec: re_ub}. Similarly, if a sample is redundant for upper bounding in a node, it will remain redundant in all child nodes of that node, and can be eliminated from the bounds tightening calculations in the current node and its child nodes. In this way, sample reduction can not only delete redundant samples at every $i_{sr}$ iterations, but also eliminate redundant information in the current node and its child nodes, thereby accelerating the overall calculation.

\subsection{Parallelization} \label{sec: parallel}
\begin{figure}[tbh]
\begin{center}
\centerline{\includegraphics[width=6in]{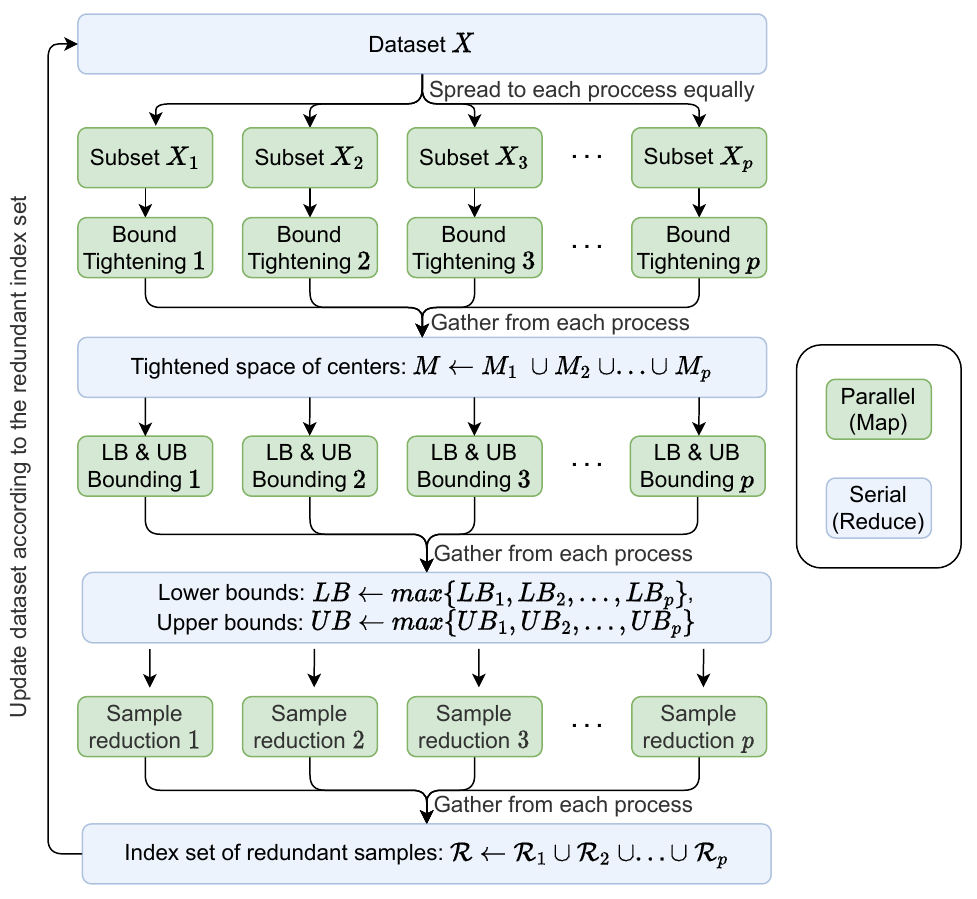}}
\caption{Parallelization of the reduced-space branch and bound scheme}
\label{fig:parallel}
\end{center}
\vskip -0.3in
\end{figure}

We also provide a parallel implementation of the whole algorithm to accelerate the solving process. Since our algorithm is primarily executed at the sample level, like $\beta_s(M)$ in the lower bounding, we can parallelize the algorithm by distributing the dataset to each process equally, then calculating on each process with the local dataset and communicating the results as needed. The detailed parallelization framework is shown in Figure \ref{fig:parallel}. Here, the green modules represent the parallel operations at each process, and the blue modules represent serial reduction operations. This parallelization framework is realized utilizing Message-Passing Interface (MPI) and MPI.jl by \citep{byrne_mpijl_2021}.

% \subsection{Effects on Computation}
% With a reduced search space, these bounds tightening techniques can significantly reduce the number of BB nodes to be explored. Moreover, they can also accelerate the calculation of bounds. For example, for the Lagrangian Problem \ref{eqn: node_reLD}, we need to calculate the contributions, $\rho_j$, of samples belonging to the overall feasible set of centers, $\hat{S}^+ = \hat{S}^{1+} \cup \hat{S}^{2+} \cup ... \cup \hat{S}^{K+}$. A reduced search space $\hat{M}$ means a reduced overall feasible selection set, which further leads to fewer calculations and faster speeds.

\section{Convergence Analysis} \label{sec:convergence}

As stated in Theorem \ref{theorem: conv_finite}, the branch-and-bound scheme for the $K$-center problem in Algorithm \ref{alg: bb_sche} converges to the global optimal solution after a finite step. In this section, we present the proof of this theorem. 

Specifically, the branch-and-bound scheme in Algorithm \ref{alg: bb_sche} branches on the region of centers, $\mu$, and generates a rooted tree with the search space $M_0$ at the root node. For the child node at $q$th level and $l_q$th iteration, we denote the search space as $M_{l_q}$. The search space of its child node is denoted as $M_{l_{q+1}}$ satisfying $M_{l_{q+1}} \subset M_{l_q}$. We denote the decreasing sequence from the root node with $M_0$ to the child node with $M_{l_q}$ as $\{M_{l_q}\}$. The search space of $k$th cluster center at $M_{l_q}$ is denoted as $M^k_{l_q}$. Along the branch-and-bound process, we can obtain a monotonically non-ascending upper bound sequence $\{\alpha_i\}$ and a monotonically non-descending lower bound sequence $\{\beta_i\}$. 

In the following convergence analysis, we adapt the fundamental conclusions from \citep{horst_global_2013} to our algorithm. It should be noted that the convergence of the $K$-center problem here is stronger than the convergence analysis in \citep{cao_scalable_2019} for two-stage nonlinear optimization problems or the convergence proof in \citep{hua_scalable_2021} for $K$-means clustering problem. Both \citet{cao_scalable_2019} and \citet{hua_scalable_2021} guarantee the convergence in the sense of $\lim \limits_{i\rightarrow\infty}\alpha_i = \lim \limits_{i\rightarrow\infty}\beta_i = z$. They can only produce a global $\epsilon$-optimal solution in a finite number of steps. While for the $K$-center problem, the algorithm can obtain an exact optimal solution (e.g., $\epsilon=0$) in a finite number of steps.

\begin{definition}
(Definition IV.3 \citep{horst_global_2013}) A bounding operation is called \textbf{finitely consistent} if, at every step, any unfathomed partition element can be further refined and if any decreasing sequence $\{M_{l_q}\}$ successively refined partition elements is finite.
\end{definition}

\begin{lemma} \label{finitely}
The bounding operation in Algorithm \ref{alg: bb_sche} is finitely consistent.
\end{lemma}
\textit{Proof.} Firstly, we prove that any unfathomed partition element $M_{l_q}$ can be further refined. Any unfathomed $M_{l_q}$ satisfies two conditions: (1) $\exists |X\cap M^k_{l_q}| > 1, k\in\mathcal{K}$, and (2) $\alpha_l - \beta(M_{l_q}) > \epsilon, \epsilon>0$. Obviously, there exists at least one partition to be further refined.

We then prove any decreasing sequences $\{M_{l_q}\}$ successively refined partition elements are finite. Assuming by contradiction that a sequence $\{M_{l_q}\}$ is infinite. In our algorithm, since we branch on the first-stage variable $\mu$ corresponding to the diameter of $M$, this subdivision is exhaustive. Therefore, we have $\lim\limits_{q\to \infty} \delta(M_{l_q})= 0$ and $\{M_{l_q}\}$ converge to one point $\bar{\mu}$ at each cluster, where $\delta(M_{l_q})$ is the the diameter of set $M_{l_q}$. 

If this point $\bar{\mu}\in X$, there exists a ball around $\bar{\mu}$, denoted as $B_r(\bar{\mu})=\{\mu\ |\ ||\mu-\bar{\mu}|| \leq r\}$, fulfilling $|X\cap B_r(\bar{\mu})| = 1$. There exists a level $q_0$ that $M_{l_q}\subset B_r(\bar{\mu}), \forall q \geq q_0$. At this $l_{q_0}$th iteration, according to the terminal conditions  $|X\cap M^k_{l_q}| = 1, \forall k \in \mathcal{K}$, the partition elements $M_{l_{q_0}}$ will not be branched anymore. Because the dataset $X$ is finite, we have the sequence $\{M_{l_q}\}$ is finite in this case. If $\bar{\mu}\not\subset X$, there is a ball around $\bar{\mu}$, denoted as $B_r(\bar{\mu})=\{\mu\ |\ ||\mu-\bar{\mu}|| \leq r\}$, satisfying $|X\cap B_r(\bar{\mu})| = 0$. There exists a level $q_0$ that $M_{l_q}\subset B_r(\bar{\mu}), \forall q \geq q_0$. At this $l_{q_0}$th iteration, $M_{l_{q_0}}$ will be deleted according to the terminal conditions. Consequently, the sequence $\{M_{l_q}\}$ is also finite in this case. In conclusion, it is impossible to exist a sequence $\{M_{l_q}\}$ that is infinite. 

\begin{theorem} \label{terminate}
(Theorem IV.1 \citep{horst_global_2013}) In a BB procedure, suppose that the bounding operation is finitely consistent. Then the procedure terminates after finitely many steps.
\end{theorem}

\begin{lemma} \label{terminate_lemma}
Algorithm \ref{alg: bb_sche} terminates after finitely many steps.
\end{lemma}
\textit{Proof.} From Lemma \ref{finitely}, the bounding operation in Algorithm \ref{alg: bb_sche} is finitely consistent. According to Theorem \ref{terminate}, we have Algorithm \ref{alg: bb_sche} terminates after finitely many steps

Finally, we prove that the BB scheme for the $K$-center problem is convergent:

\noindent\textbf{Theorem \ref{theorem: conv_finite}. }\textit{Algorithm \ref{alg: bb_sche} is convergent to the global optimal solution after a finite step $L$, with $\beta_L=z=\alpha_L$, by only branching on the space of $\mu$.}

\textit{Proof.} From Lemma \ref{terminate_lemma}, Algorithm \ref{alg: bb_sche} terminates after finite steps. The algorithm terminates with two situations. The first situations is $|\beta_l - \alpha_l| \leq \epsilon, \epsilon \geq 0$. When $\epsilon$ is set to be 0, we have $\beta_l=z=\alpha_l$.

The second situation is the branch node set $\mathbb{M} = \emptyset$. A branch node with $M$ is deleted from $\mathbb{M}$ and not further partitioned if it satisfies  $\beta(M) > \alpha_l$ or $|X\cap M^k| = 1, \forall k\in \mathcal{K}$. In the first case, it is obvious that this branch node does not contain the global optimal solution $\mu^*$. Therefore, the branch node with $M'$ containing the optimal solution $\mu^*$ is not further partitioned because the second case $|X\cap M'^k| = 1, \forall k\in \mathcal{K}$. After bounds tightening according to the ``centers on samples'' constraint, the tightened node $M'= \{\mu^*\}$. Obviously for this tightened node, we have $\beta_l =\beta(M') = z = \alpha(M')= \alpha_l$.  In this way, we have proved Theorem \ref{theorem: conv_finite}.

\section{Numerical Results} \label{sec:results}
In this section, we report the detailed implementation of our algorithm and the numerical results on synthetic and real-world datasets.

\subsection{Implementation Details}
We denote our tailored reduced-space branch and bound algorithm \ref{alg: bb_sche} with and without acceleration techniques as \texttt{BB+CF+BT} and \texttt{BB+CF} correspondingly. All our algorithms are implemented in \texttt{Julia}, and the parallel version is realized using Message Passing Interface through the \texttt{MPI.jl} module. We compare the performance of our algorithm with the state-of-art global optimizer \texttt{CPLEX 20.1.0} \citep{cplex_v2010_2020} and the heuristic algorithm, Farthest First Traversal (\texttt{FFT}) as shown in Algorithm \ref{alg: fft}. The initial points severely influence the results of \texttt{FFT}. Therefore, we execute \texttt{FFT} for 100 trails with randomly selected initial points and report the best results. As for CPLEX, we use the MINLP formulation \ref{eqn:overall} with the symmetry-breaking constraints to solve the $K$-center problem. We executed all experiments on the high-performance computing cluster Niagara in the Digital Research Alliance of Canada. Each computing node of the Niagara cluster has 40 Intel ``Skylake'' cores and 188 GiB of RAM. For the global optimizer \texttt{CPLEX} and our algorithms, a time limit of 4 hours is set to compare the performance fairly and avoid unacceptable computational costs. For our algorithms, there is also an optimality gap limit of $0.1\%$. The source code is available as Supplemental Material at \url{https://doi.org/10.1287/mnsc.2023.00218}.

In order to evaluate the performance extensively, we execute all the algorithms on both synthetic and real-world datasets. The synthetic datasets are generated using \texttt{Distributions.jl} and \texttt{Random.jl} modules in \texttt{Julia}. We generate the synthetic datasets with 3 Gaussian clusters, 2 attributes, and varying numbers of samples. As for the real-world datasets, we use 30 datasets from the UCI Machine Learning Repository \citep{Dua_2017}, datasets \texttt{Pr2392} from \citep{padberg_branch-and-cut_1991}, \texttt{Hemi} from \citep{MLbio} and \texttt{Taxi} from \citep{schneider_analyzing_2015}. The number of samples ranges from 150 to 1,120,841,769. The number of attributes ranges from 2 to 68. The detailed characteristics of datasets can be found in the following result tables. 

We report four criteria in the following result tables to compare the performance of algorithms: upper bound (\texttt{UB}), optimality gap (\texttt{Gap}), the number of solved BB nodes (\texttt{Nodes}), and the run time (\texttt{Time}). \texttt{UB} is the best objective value of the $K$-center Problem \ref{eq:obj}. \texttt{Gap} represents the relative difference between the best lower bound (LB) and UB. It is defined as $\texttt{Gap} = \frac{UB-LB}{UB}
\times 100\%$. The optimality gap is a unique property of the deterministic global optimization algorithm. The heuristic algorithm (\texttt{FFT}) does not have this property. \texttt{Nodes} and \texttt{Time} are the iteration number and the run time of the BB scheme from the beginning to the termination. 

\subsection{Serial Results on Synthetic Datasets}\label{sec:result_large}

Table \ref{tab:syn_data} reports the serial results of synthetic datasets with different numbers of samples and different desired clusters ($K=3, 5, 10$). Compared with the heuristic method \texttt{FFT}, our algorithm \texttt{BB+LD+BT} can reduce \texttt{UB} by 29.4\% average on these synthetic datasets. These results validate the conclusion from \cite{garcia-diaz_approximation_2019} that these 2-approximation heuristic algorithms perform poorly in practice despite the solution quality guarantee. 

As for the comparison of global optimizers, the direct usage of \texttt{CPLEX} on Problem \ref{eqn:overall} could not converge to a small optimality gap$\leq 0.1\%$ within 4 hours on all the synthetic datasets. \texttt{BB+LD} without acceleration techniques can obtain the small optimality gap$\leq 0.1\%$ within 4 hours on  synthetic datasets smaller than 42,000 samples with desired clusters $K=3$. The algorithm \texttt{BB+LD+BT} can obtain the best upper bounds and reach a satisfactory gap$\leq 0.1\%$ in most experiments within 4 hours. Moreover, compared with \texttt{BB+LD}, \texttt{BB+LD+BT} needs fewer nodes and less run time to obtain the same optimality gap. For example, for the \texttt{Syn-1200} dataset with $K=3$, \texttt{BB+LD} need 1,155,375 nodes and 3609 seconds to reach a gap$\leq0.1\%$, while \texttt{BB+LD+BT} only needs 23 nodes and 13.5 seconds. These comparisons between \texttt{BB+LD} and \texttt{BB+LD+BT} demonstrate the acceleration techniques in Section \ref{sec:acc} can significantly reduce the search space and accelerate the BB procedure.

\subsection{Serial Results on Real-world datasets}\label{sec:result_small}

Table \ref{tab:s_results}, Table \ref{tab:l_results}, and Table \ref{tab:mil_result} show the serial results on real-world datasets with different sample numbers and desired cluster numbers ($K=3,5,10$). In these tables, we highlight the best results among these algorithms with the optimality gap$\leq0.1\%$. These real-world results are consistent with the results of synthetic datasets. 

The best solutions generated by the heuristic method (\texttt{FFT}) can be far from optimal in these tables, even for very small datasets. For example, for \texttt{IRIS} dataset, \texttt{FFT} obtains \texttt{UB} of 3.66 while our algorithm and \texttt{CPLEX} give a \texttt{UB} of 2.04 with $\leq 0.1\%$ gap. Compared with \texttt{FFT}, our algorithm \texttt{BB+CF+BT} can averagely reduce the \texttt{UB} by $22.2\%$ on these real-world datasets and $25.8\%$ on all the synthetic and real-world datasets. Even for experiments terminated with large gaps, in most cases, \texttt{BB+CF+BT} can obtain a smaller \texttt{UB} than \texttt{FFT}.

For small datasets, our algorithms \texttt{BB+CF} and \texttt{BB+CF+BT} can obtain the same \texttt{UB} as \texttt{CPLEX}. However, \texttt{CPLEX} needs significantly more run time and nodes than our algorithms. For all datasets with more than 740 samples, \texttt{CPLEX} cannot even give an optimality gap$\leq 50\%$ within 4 hours. On the contrary, \texttt{BB+CF+BT} can obtain the best \texttt{UB} and a satisfactory gap$\leq 0.1\%$ for most datasets. 

The comparisons of the two versions of our algorithms \texttt{BB+CF} and \texttt{BB+CF+BT} demonstrate that the acceleration techniques in Section \ref{sec:acc} can significantly reduce the computational time and the number of BB nodes to solve the problems. Remarkably, with these acceleration techniques, we can even solve several datasets in the root node (\texttt{Nodes}${=}$1), e.g., the datasets \texttt{iris}, \texttt{HF}, and \texttt{SGC}. Besides, \texttt{BB+CF+BT} results with superscript $^1$ in these tables mean we can assign $K$ farthest initial seeds through FFT at the root node as described in Section \ref{sec:fbbt_fft}. We can obtain the initial seeds for about half of the datasets when $K=3$. Moreover, the number of nodes is much smaller for the datasets with initial seeds than the datasets without initial seeds. This phenomenon indicates the initial seeds are essential for cluster assignment and bounds tightening since we need at least one assigned sample at each cluster to execute the sample-based assignment.

For most of the datasets with millions of samples and $K=3$ in Table \ref{tab:mil_result}, \texttt{BB+CF+BT} can converge to a small gap$\leq 0.1\%$ and provide the best optimal solution after 4 hours of running. To the best of our knowledge, it is the first time that the $K$-center problem is solved under a relatively small gap$\leq 0.1\%$ within 4 hours on datasets over \textbf{14 million} samples in the \textbf{serial mode}. 

As a drawback, our algorithm \texttt{BB+LD+BT} still struggles to obtain a small optimality gap when the desired number of clusters is larger than 3. However, it should be noted the state-of-art global optimizer \texttt{CPLEX} cannot even solve any datasets to gap$\leq50\%$ when $K>3$. On the contrary, \texttt{BB+LD+BT} can obtain gap$\leq0.1\%$ on most datasets with less than 5 million samples and $K=5$. Moreover, for the cases when our algorithm \texttt{BB+LD+BT} cannot obtain a small optimality gap, it still gives the best \texttt{UB} among all the algorithms in these experiments.

\subsection{Parallel Results on Huge-scale Real-world datasets}\label{sec:result_parallel}
To fully utilize the computational ability of high-performance clusters, we implement our algorithm \texttt{BB+CF+BT} in a parallel manner as shown in Section \ref{sec: parallel}. Here, we test the parallel algorithm on datasets that couldn't obtain a small gap$\leq0.1\%$ for $K=3$ within 4 hours in the serial mode, including two datasets with ten million samples, \texttt{HIGGS} and \texttt{BigCross}. Moreover, we also extend the experiments to a billion-scale dataset called \texttt{Taxi}. This billion-scale dataset contains over 1.1 billion individual taxi trips with 12 attributes in New York City from January 2009 through June 2015. We preprocess the Taxi dataset according to the analysis by \cite{schneider_analyzing_2015} to remove outliers and missing values in the dataset. As an outcome shown in Table \ref{tbl: para_result}, the parallel version of \texttt{BB+CF+BT} can reach a small optimality gap$\leq 0.1\%$ and a better \texttt{UB} on the datasets \texttt{BigCross} and \texttt{Taxi} within 4 hours. For the dataset \texttt{HIGGS}, the parallel version achieves a smaller \texttt{UB} and gap compared to the heuristic method and the serial version. As far as we know, this is the first time that the $K$-center problem is solved under a relatively small gap$\leq 0.1\%$ within 4 hours on the \textbf{billion-scale dataset}.

\section{Conclusion}
We propose a global optimization algorithm for the $K$-center problem using a tailored reduced space branch and bound scheme. In this algorithm, we only need to branch on the region of cluster centers to guarantee convergence to the global optimal solution in a finite step. 

We give a two-stage decomposable formulation and an MINLP formulation of the $K$-center problem. With this two-stage formulation, we develop a lower bound with closed-form solutions by relaxing the non-anticipativity constraints and the ``centers on sample'' constraints. As an outcome, the proposed bounding methods are extremely computationally efficient with no needs to solve any optimization sub-problems using any optimizers.

Along with the BB procedure, we introduce several acceleration techniques based on the MINLP formulation, including bounds tightening, and sample reduction. Numerical experiments show these acceleration techniques can significantly reduce the search space and accelerate the solving procedure. Moreover, we also give a parallel implementation of our algorithm to fully utilize the computational power of modern high performance clusters.

Extensive numerical experiments have been conducted on synthetic and real-world datasets. These results exhibit the efficiency of our algorithm: we can solve the real-world datasets with up to \textbf{ten million samples} in the serial mode and \textbf{one billion samples} in the parallel mode to a small optimality gap ($\leq$0.1\%) within 4 hours.

Finally, we also declare that our algorithm is promised to extend to deal with certain constrained versions of $K$-center problems. For example, the capacitated restricted version, absolute and vertex restricted version \citep{calik_exact_2013}. We are interested in developing these variants in future work. 

% Acknowledgments here

% References here (outcomment the appropriate case)
\bibliographystyle{informs2014} 
\bibliography{references}

\newpage

\begin{table}[!htp]
\caption{Serial results on synthetic datasets}
\tiny
\label{tab:syn_data}
\setlength\tabcolsep{2pt}
\centering
\begin{threeparttable}
\renewcommand{\arraystretch}{1.5}
\begin{tabular}{cccccccccccccccc}
\hline
\multirow{2}{*}{Dataset} & \multirow{2}{*}{\begin{tabular}[c]{@{}c@{}}Sam\\ ple\end{tabular}} & \multirow{2}{*}{\begin{tabular}[c]{@{}c@{}}Dimen\\ sion\end{tabular}} & \multirow{2}{*}{Method} & \multicolumn{4}{c}{K=3} & \multicolumn{4}{c}{K=5} & \multicolumn{4}{c}{K=10} \\ \cline{5-16} 
 &  &  &  & UB & Nodes & \begin{tabular}[c]{@{}c@{}}Gap\\      (\%)\end{tabular} & \begin{tabular}[c]{@{}c@{}}Time\\      (s)\end{tabular} & UB & Nodes & \begin{tabular}[c]{@{}c@{}}Gap\\      (\%)\end{tabular} & \begin{tabular}[c]{@{}c@{}}Time\\      (s)\end{tabular} & UB & Nodes & \begin{tabular}[c]{@{}c@{}}Gap\\      (\%)\end{tabular} & \begin{tabular}[c]{@{}c@{}}Time\\      (s)\end{tabular} \\
 \hline
\multirow{4}{*}{Syn-300} & \multirow{4}{*}{3.0E+2} & \multirow{4}{*}{2} & FFT & 69.68 & - & - & - & 43.33 & - & - & - & 21.88 & - & - & - \\
 &  &  & CPLEX & 61.75 & 2.9E+4 & $\leq$0.1 & 29 & 37.14 & 2.3E+7 & 19.4 & 4h & 16.06 & 1.2E+7 & 100.0 & 4h \\
 &  &  & BB+CF & 61.75 & 5.5E+4 & $\leq$0.1 & 46 & 37.14 & 2.3E+6 & 16.2 & 4h & 15.64 & 1.7E+6 & 100.0 & 4h \\
 &  &  & \textbf{BB+CF+BT} & \textbf{61.75} & \textbf{17} & \textbf{$\leq$0.1} & \textbf{13} & \textbf{37.14} & \textbf{1,764} & \textbf{$\leq$0.1} & \textbf{15} & \textbf{12.31} & \textbf{2.0E+4} & \textbf{$\leq$0.1} & \textbf{38} \\
 \hline
\multirow{4}{*}{Syn-1200} & \multirow{4}{*}{1.2E+3} & \multirow{4}{*}{2} & FFT & 93.34 & - & - & - & 58.46 & - & - & - & 30.49 & - & - & - \\
 &  &  & CPLEX & 84.81 & 5.8E+6 & 1.6 & 4h & 34.29 & 3.5E+6 & 7.8 & 4h & 89.32 & 8.1E+5 & 100.0 & 4h \\
 &  &  & BB+CF & 84.81 & 1.2E+6 & $\leq$0.1 & 3,609 & 34.29 & 1.4E+6 & 12.5 & 4h & 21.81 & 1.0E+6 & 100.0 & 4h \\
 &  &  & \textbf{BB+CF+BT} & \textbf{84.81} & \textbf{23} & \textbf{$\leq$0.1}$^1$ & \textbf{14} & \textbf{34.29} & \textbf{411} & \textbf{$\leq$0.1} & \textbf{15} & \textbf{14.51} & \textbf{3.0E+4} & \textbf{$\leq$0.1} & \textbf{148} \\
 \hline
\multirow{4}{*}{Syn-2100} & \multirow{4}{*}{2.1E+3} & \multirow{4}{*}{2} & FFT & 106.50 & - & - & - & 72.70 & - & - & - & 36.04 & - & - & - \\
 &  &  & CPLEX & 95.10 & 3.0E+6 & 0.2 & 4h & 49.32 & 1.3E+6 & 100.0 & 4h & 193.26 & 3.4E+5 & 100.0 & 4h \\
 &  &  & BB+CF & 95.10 & 1.5E+6 & $\leq$0.1 & 11,606 & 42.58 & 1.0E+6 & 20.8 & 4h & 25.78 & 5.3E+5 & 100.0 & 4h \\
 &  &  & \textbf{BB+CF+BT} & \textbf{95.10} & \textbf{17} & \textbf{$\leq$0.1}$^1$ & \textbf{13} & \textbf{42.58} & \textbf{455} & \textbf{$\leq$0.1} & \textbf{16} & \textbf{17.65} & \textbf{8.9E+4} & \textbf{$\leq$0.1} & \textbf{725} \\
 \hline
\multirow{4}{*}{Syn-42000} & \multirow{4}{*}{4.2E+4} & \multirow{4}{*}{2} & FFT & 161.98 & - & - & - & 96.12 & - & - & - & 47.21 & - & - & - \\
 &  &  & CPLEX & \multicolumn{4}{c}{No feasible solution} & \multicolumn{4}{c}{No feasible solution} & \multicolumn{4}{c}{No feasible solution} \\
 &  &  & BB+CF & 142.33 & 1.7E+5 & 6.7 & 4h & 63.40 & 1.0E+5 & 28.1 & 4h & 44.24 & 5.4E+4 & 100.0 & 4h \\
 &  &  & \textbf{BB+CF+BT} & \textbf{142.33} & \textbf{103} & \textbf{$\leq$0.1} & \textbf{21} & \textbf{62.77} & \textbf{5.0E+3} & \textbf{$\leq$0.1} & \textbf{363} & 28.29 & 5.8E+4 & 36.1 & 4h \\
 \hline
\multirow{4}{*}{Syn-210000} & \multirow{4}{*}{2.1E+5} & \multirow{4}{*}{2} & FFT & 175.81 & - & - & - & 120.78 & - & - & - & 66.79 & - & - & - \\
 &  &  & CPLEX & \multicolumn{4}{c}{No feasible solution} & \multicolumn{4}{c}{No feasible solution} & \multicolumn{4}{c}{No feasible solution} \\
 &  &  & BB+CF & 168.57 & 4.4E+4 & 7.0 & 4h & 77.02 & 2.5E+4 & 43.8 & 4h & 53.73 & 1.4E+4 & 100.0 & 4h \\
 &  &  & \textbf{BB+CF+BT} & \textbf{168.57} & \textbf{5} & \textbf{$\leq$0.1}$^1$ & \textbf{21} & \textbf{71.88} & \textbf{2.4E+3} & \textbf{$\leq$0.1} & \textbf{1,118} & 44.48 & 1.2E+4 & 72.2 & 4h\\
 \hline
\end{tabular}
\begin{tablenotes}
    \item[1] Can assign $K$ initial seeds through FFT at the root node. \texttt{BB+CF+BT} results without this superscript means can not assign initial seeds.
\end{tablenotes}
\end{threeparttable}
\end{table}

\begin{table}[!htp]
\caption{Serial results on small-scale datasets ($S\leq 1,000$)}\label{tab:s_results}
\tiny
\setlength\tabcolsep{2pt}
\begin{threeparttable}
\renewcommand{\arraystretch}{1.5}
\begin{tabular}{cccccccccccccccc}
\hline
\multirow{2}{*}{Dataset} & \multirow{2}{*}{\begin{tabular}[c]{@{}c@{}}Sam\\      ple\end{tabular}} & \multirow{2}{*}{\begin{tabular}[c]{@{}c@{}}Dimen\\      sion\end{tabular}} & \multirow{2}{*}{Method} & \multicolumn{4}{c}{K=3} & \multicolumn{4}{c}{K=5} & \multicolumn{4}{c}{K=10} \\ \cline{5-16} 
 &  &  &  & UB & Nodes & \begin{tabular}[c]{@{}c@{}}Gap\\      (\%)\end{tabular} & \begin{tabular}[c]{@{}c@{}}Time\\      (s)\end{tabular} & UB & Nodes & \begin{tabular}[c]{@{}c@{}}Gap\\      (\%)\end{tabular} & \begin{tabular}[c]{@{}c@{}}Time\\      (s)\end{tabular} & UB & Nodes & \begin{tabular}[c]{@{}c@{}}Gap\\      (\%)\end{tabular} & \begin{tabular}[c]{@{}c@{}}Time\\      (s)\end{tabular} \\
 \hline
\multirow{4}{*}{iris} & \multirow{4}{*}{150} & \multirow{4}{*}{4} & FFT & 2.65 & - & - & - & 1.80 & - & - & - & 0.95 & - & - & - \\
 &  &  & CPLEX & 2.04 & 1.2E+5 & $\leq$0.1 & 46 & 1.54 & 2.8E+6 & 60.0 & 4h & 1.21 & 1.4E+7 & 100.0 & 4h \\
 &  &  & BB+CF & 2.04 & 1.3E+4 & $\leq$0.1 & 17 & 1.20 & 3.1E+6 & $\leq$0.1 & 5,472 & 0.74 & 2.2E+6 & 100.0 & 4h \\
 &  &  & \textbf{BB+CF+BT} & \textbf{2.04} & \textbf{1} & \textbf{$\leq$0.1}$^1$ & \textbf{12} & \textbf{1.20} & \textbf{409} & \textbf{$\leq$0.1} & \textbf{14} & 0.66 & 9.6E+5 & 25.8 & 4h \\
  \hline
\multirow{4}{*}{seeds} & \multirow{4}{*}{210} & \multirow{4}{*}{7} & FFT & 13.17 & - & - & - & 9.01 & - & - & - & 4.48 & - & - & - \\
 &  &  & CPLEX & 10.44 & 1.2E+6 & $\leq$0.1 & 542 & 11.61 & 2.5E+6 & 96.1 & 4h & 21.48 & 5.6E+5 & 100.0 & 4h \\
 &  &  & BB+CF & 10.44 & 7.2E+3 & $\leq$0.1 & 17 & 7.22 & 2.7E+6 & 8.3 & 4h & 3.51 & 1.5E+6 & 100.0 & 4h \\
 &  &  & \textbf{BB+CF+BT} & \textbf{10.44} & \textbf{21} & \textbf{$\leq$0.1}$^1$ & \textbf{13} & \textbf{7.22} & \textbf{1,444} & \textbf{$\leq$0.1} & \textbf{15} & \textbf{2.92} & \textbf{2.1E+5} & \textbf{$\leq$0.1} & \textbf{569} \\
  \hline
\multirow{4}{*}{glass} & \multirow{4}{*}{214} & \multirow{4}{*}{9} & FFT & 27.52 & - & - & - & 22.28 & - & - & - & 11.73 & - & - & - \\
 &  &  & CPLEX & \multicolumn{4}{c}{Out of memory} & \multicolumn{4}{c}{Out of memory} & \multicolumn{4}{c}{Out of memory} \\
 &  &  & BB+CF & 27.52 & 5.6E+3 & $\leq$0.1 & 15 & 16.44 & 9.7E+5 & $\leq$0.1 & 1,522 & 10.64 & 1.4E+6 & 100.0 & 4h \\
 &  &  & \textbf{BB+CF+BT} & \textbf{27.52} & \textbf{191} & \textbf{$\leq$0.1} & \textbf{13} & \textbf{16.44} & \textbf{4.4E+3} & \textbf{$\leq$0.1} & \textbf{17} & \textbf{7.95} & \textbf{1.7E+6} & \textbf{$\leq$0.1} & \textbf{9,180} \\
  \hline
\multirow{4}{*}{BM} & \multirow{4}{*}{249} & \multirow{4}{*}{6} & FFT & 1.52E+04 & - & - & - & 1.12E+04 & - & - & - & 5.33E+03 & - & - & - \\
 &  &  & CPLEX & \multicolumn{4}{c}{No feasible solution} & 1.48E+04 & 8.8E+6 & 100.0 & 4h & 1.63E+04 & 2.4E+6 & 100.0 & 4h \\
 &  &  & BB+CF & 1.05E+04 & 1.4E+4 & $\leq$0.1 & 22 & 6.32E+03 & 2.2E+6 & 12.0 & 4h & 5.01E+03 & 1.4E+6 & 100.0 & 4h \\
 &  &  & \textbf{BB+CF+BT} & \textbf{1.05E+04} & \textbf{63} & \textbf{$\leq$0.1}$^1$ & \textbf{13} & \textbf{6.32E+03} & \textbf{1.8E+4} & \textbf{$\leq$0.1} & \textbf{29} & 4.98E+03 & 6.7E+5 & 97.9 & 4h \\
  \hline
\multirow{4}{*}{UK} & \multirow{4}{*}{258} & \multirow{4}{*}{5} & FFT & 0.70 & - & - & - & 0.57 & - & - & - & 0.42 & - & - & - \\
 &  &  & CPLEX & \multicolumn{4}{c}{Out of memory} & \multicolumn{4}{c}{Out of memory} & \multicolumn{4}{c}{Out of memory} \\
 &  &  & BB+CF & 0.53 & 3.2E+5 & $\leq$0.1 & 258 & 0.43 & 1.5E+6 & 43.9 & 4h & 0.33 & 1.4E+6 & 100.0 & 4h \\
 &  &  & \textbf{BB+CF+BT} & \textbf{0.53} & \textbf{1.6E+4} & \textbf{$\leq$0.1} & \textbf{23} & 0.43 & 8.9E+5 & 26.9 & 4h & 0.31 & 6.1E+5 & 97.3 & 4h \\
  \hline
\multirow{4}{*}{HF} & \multirow{4}{*}{299} & \multirow{4}{*}{12} & FFT & 2.69E+10 & - & - & - & 1.17E+10 & - & - & - & 1.68E+09 & - & - & - \\
 &  &  & CPLEX & \multicolumn{4}{c}{No feasible solution} & \multicolumn{4}{c}{No feasible solution} & \multicolumn{4}{c}{No feasible solution} \\
 &  &  & BB+CF & 1.72E+10 & 339 & $\leq$0.1 & 10 & 1.02E+10 & 2.1E+4 & $\leq$0.1 & 44 & 1.52E+09 & 3.4E+6 & 100.0 & 4h \\
 &  &  & \textbf{BB+CF+BT} & \textbf{1.72E+10} & \textbf{1} & \textbf{$\leq$0.1}$^1$ & \textbf{12} & \textbf{1.02E+10} & \textbf{557} & \textbf{$\leq$0.1} & \textbf{14} & 1.44E+09 & 1.2E+6 & 53.2 & 4h \\
  \hline
\multirow{4}{*}{Who} & \multirow{4}{*}{440} & \multirow{4}{*}{8} & FFT & 4.58E+09 & - & - & - & 3.18E+09 & - & - & - & 9.81E+08 & - & - & - \\
 &  &  & CPLEX & \multicolumn{4}{c}{No feasible solution} & \multicolumn{4}{c}{No feasible solution} & \multicolumn{4}{c}{No feasible solution} \\
 &  &  & BB+CF & 3.49E+09 & 3.4E+3 & $\leq$0.1 & 15 & 2.11E+09 & 1.7E+5 & $\leq$0.1 & 341 & 9.27E+08 & 1.5E+6 & 100.0 & 4h \\
 &  &  & \textbf{BB+CF+BT} & \textbf{3.49E+09} & \textbf{375} & \textbf{$\leq$0.1} & \textbf{14} & \textbf{2.11E+09} & \textbf{2.3E+3} & \textbf{$\leq$0.1} & \textbf{16} & 8.21E+08 & 8.4E+5 & 62.0 & 4h \\
  \hline
\multirow{4}{*}{HCV} & \multirow{4}{*}{602} & \multirow{4}{*}{12} & FFT & 1.75E+05 & - & - & - & 8.38E+04 & - & - & - & 3.03E+04 & - & - & - \\
 &  &  & CPLEX & 1.41E+05 & 9.5E+5 & $\leq$0.1 & 3,720 & 8.73E+04 & 5.1E+5 & 100.0 & 4h & 4.47E+04 & 3.8E+5 & 100.0 & 4h \\
 &  &  & BB+CF & 1.41E+05 & 291 & $\leq$0.1 & 10 & 6.37E+04 & 2.2E+4 & $\leq$0.1 & 76 & 2.36E+04 & 1.4E+6 & 100.0 & 4h \\
 &  &  & \textbf{BB+CF+BT} & \textbf{1.41E+05} & \textbf{39} & \textbf{$\leq$0.1} & \textbf{13} & \textbf{6.37E+04} & \textbf{583} & \textbf{$\leq$0.1} & \textbf{15} & \textbf{2.16E+04} & \textbf{7.6E+5} & \textbf{$\leq$0.1} & \textbf{6,300} \\
  \hline
\multirow{4}{*}{Abs} & \multirow{4}{*}{740} & \multirow{4}{*}{21} & FFT & 1.94E+04 & - & - & - & 1.19E+04 & - & - & - & 7.81E+03 & - & - & - \\
 &  &  & CPLEX & 1.72E+04 & 9.3E+5 & 52.0 & 4h & 3.02E+04 & 1.7E+4 & 100.0 & 4h & \multicolumn{4}{c}{No feasible solution} \\
 &  &  & BB+CF & 1.39E+04 & 3.3E+4 & $\leq$0.1 & 153 & 9.93E+03 & 1.5E+6 & 18.9 & 4h & 6.15E+03 & 9.8E+5 & 100.0 & 4h \\
 &  &  & \textbf{BB+CF+BT} & \textbf{1.39E+04} & \textbf{611} & \textbf{$\leq$0.1} & \textbf{15} & \textbf{9.92E+03} & \textbf{5.0E+4} & \textbf{$\leq$0.1} & \textbf{178} & 6.37E+03 & 4.6E+5 & 98.3 & 4h \\
  \hline
\multirow{4}{*}{TR} & \multirow{4}{*}{980} & \multirow{4}{*}{10} & FFT & 7.32 & - & - & - & 6.42 & - & - & - & 4.41 & - & - & - \\
 &  &  & CPLEX & 8.32 & 1.6E+6 & 54.5 & 4h & 7.82 & 2.0E+5 & 100.0 & 4h & 8.70 & 3.3E+4 & 100.0 & 4h \\
 &  &  & BB+CF & 5.94 & 7.4E+5 & $\leq$0.1 & 2,953 & 4.49 & 1.3E+6 & 47.7 & 4h & 3.69 & 9.8E+5 & 100.0 & 4h \\
 &  &  & \textbf{BB+CF+BT} & \textbf{5.94} & \textbf{3.3E+4} & \textbf{$\leq$0.1} & \textbf{83} & 4.49 & 1.0E+6 & 24.6 & 4h & 3.73 & 5.0E+5 & 99.9 & 4h \\
  \hline
\multirow{4}{*}{SGC} & \multirow{4}{*}{1,000} & \multirow{4}{*}{21} & FFT & 1.33E+07 & - & - & - & 4.08E+06 & - & - & - & 9.50E+05 & - & - & - \\
 &  &  & CPLEX & 9.45E+06 & 5.0E+4 & 100.0 & 4h & 1.56E+08 & 10 & 100.0 & 4h & \multicolumn{4}{c}{No feasible solution} \\
 &  &  & BB+CF & 9.45E+06 & 411 & $\leq$0.1 & 12 & 3.91E+06 & 2.8E+4 & $\leq$0.1 & 185 & 9.50E+05 & 9.6E+5 & 100.0 & 4h \\
 &  &  & \textbf{BB+CF+BT} & \textbf{9.45E+06} & \textbf{1} & \textbf{$\leq$0.1}$^1$ & \textbf{12} & \textbf{3.91E+06} & \textbf{1} & \textbf{$\leq$0.1}$^1$ & \textbf{12} & 9.50E+05 & 5.8E+5 & 100.0 & 4h \\
  \hline
\end{tabular}
\begin{tablenotes}
    \item[1] Can assign $K$ initial seeds through FFT at the root node. \texttt{BB+CF+BT} results without this superscript means can not assign initial seeds.
  \end{tablenotes}
\end{threeparttable}
\end{table}

\begin{table}[!htp]
\caption{Serial results on large-scale datasets ($1,000{<}S{<}1,000,000$)}
\tiny
\label{tab:l_results}
\setlength\tabcolsep{1pt}
\begin{center}
\begin{threeparttable}
\renewcommand{\arraystretch}{1.5}
\begin{tabular}{cccccccccccccccc}
\hline
\multirow{2}{*}{Dataset} & \multirow{2}{*}{\begin{tabular}[c]{@{}c@{}}Sam\\      ple\end{tabular}} & \multirow{2}{*}{\begin{tabular}[c]{@{}c@{}}Dimen\\      sion\end{tabular}} & \multirow{2}{*}{Method} & \multicolumn{4}{c}{K=3} & \multicolumn{4}{c}{K=5} & \multicolumn{4}{c}{K=10} \\ \cline{5-16} 
 &  &  &  & UB & Nodes & \begin{tabular}[c]{@{}c@{}}Gap\\      (\%)\end{tabular} & \begin{tabular}[c]{@{}c@{}}Time\\      (s)\end{tabular} & UB & Nodes & \begin{tabular}[c]{@{}c@{}}Gap\\      (\%)\end{tabular} & \begin{tabular}[c]{@{}c@{}}Time\\      (s)\end{tabular} & UB & Nodes & \begin{tabular}[c]{@{}c@{}}Gap\\      (\%)\end{tabular} & \begin{tabular}[c]{@{}c@{}}Time\\      (s)\end{tabular} \\
 \hline
\multirow{4}{*}{hemi} & \multirow{4}{*}{1,955} & \multirow{4}{*}{7} & FFT & 1.06E+05 & - & - & - & 3.31E+04 & - & - & - & 1.42E+04 & - & - & - \\
 &  &  & CPLEX & 4.08E+05 & 1.4E+5 & 100.0 & 4h & 4.08E+05 & 1.2E+5 & 100.0 & 4h & 4.08E+05 & 5 & 100.0 & 4h \\
 &  &  & BB+CF & 4.08E+05 & 1.4E+5 & 50.0 & 4h & 2.18E+04 & 6.8E+5 & $\leq$0.1 & 4,212 & 1.42E+04 & 7.1E+5 & 100.0 & 4h \\
 &  &  & \textbf{BB+CF+BT} & \textbf{6.49E+04} & \textbf{11} & \textbf{$\leq$0.1}$^1$ & \textbf{14} & \textbf{2.18E+04} & \textbf{158} & \textbf{$\leq$0.1} & \textbf{14} & \textbf{7.20E+03} & \textbf{8.2E+3} & \textbf{$\leq$0.1} & \textbf{89} \\
 \hline
\multirow{4}{*}{pr2392} & \multirow{4}{*}{2,392} & \multirow{4}{*}{2} & FFT & 3.75E+07 & - & - & - & 2.11E+07 & - & - & - & 1.02E+07 & - & - & - \\
 &  &  & CPLEX & 6.51E+07 & 2.1E+5 & 100.0 & 4h & 5.66E+07 & 3.6E+4 & 100.0 & 4h & 4.23E+07 & 4.8E+4 & 100.0 & 4h \\
 &  &  & BB+CF & 2.93E+07 & 5.9E+4 & $\leq$0.1 & 297 & 1.52E+07 & 8.3E+5 & 20.9 & 4h & 1.02E+07 & 6.0E+5 & 100.0 & 4h \\
 &  &  & \textbf{BB+CF+BT} & \textbf{2.93E+07} & \textbf{207} & \textbf{$\leq$0.1} & \textbf{15} & \textbf{1.46E+07} & \textbf{6.6E+3} & \textbf{$\leq$0.1} & \textbf{45} & 8.70E+06 & 4.2E+5 & 59.6 & 4h \\
 \hline
\multirow{4}{*}{TRR} & \multirow{4}{*}{5,454} & \multirow{4}{*}{24} & FFT & 101.55 & - & - & - & 95.26 & - & - & - & 86.87 & - & - & - \\
 &  &  & CPLEX & 166.61 & 1 & 100.0 & 4h & \multicolumn{4}{c}{No feasible solution} & \multicolumn{4}{c}{No feasible solution} \\
 &  &  & BB+CF & 89.78 & 3.6E+5 & 73.1 & 4h & 85.07 & 2.7E+5 & 100.0 & 4h & 78.53 & 1.7E+5 & 100.0 & 4h \\
 &  &  & BB+CF+BT & 88.30 & 2.8E+5 & 64.2 & 4h & 84.80 & 2.0E+5 & 100.0 & 4h & 77.37 & 1.3E+5 & 100.0 & 4h \\
 \hline
\multirow{4}{*}{AC} & \multirow{4}{*}{7,195} & \multirow{4}{*}{22} & FFT & 3.59 & - & - & - & 2.79 & - & - & - & 2.28 & - & - & - \\
 &  &  & CPLEX & \multicolumn{4}{c}{No feasible solution} & \multicolumn{4}{c}{No feasible solution} & \multicolumn{4}{c}{No feasible solution} \\
 &  &  & BB+CF & 2.75 & 3.1E+5 & 42.6 & 4h & 2.26 & 2.3E+5 & 72.0 & 4h & 2.14 & 1.4E+5 & 100.0 & 4h \\
 &  &  & BB+CF+BT & 2.78 & 2.7E+5 & 38.9 & 4h & 2.26 & 1.8E+5 & 70.2 & 4h & 1.90 & 1.2E+5 & 100.0 & 4h \\
 \hline
\multirow{4}{*}{rds\_cnt} & \multirow{4}{*}{10,000} & \multirow{4}{*}{4} & FFT & 1.93E+04 & - & - & - & 5.93E+03 & - & - & - & 1.44E+03 & - & - & - \\
 &  &  & CPLEX & 6.86E+04 & 3.2E+4 & 100.0 & 4h & \multicolumn{4}{c}{No feasible solution} & \multicolumn{4}{c}{No feasible solution} \\
 &  &  & BB+CF & 1.39E+04 & 639 & $\leq$0.1 & 25 & 4.90E+03 & 1.6E+5 & $\leq$0.1 & 6,048 & 1.44E+03 & 1.8E+5 & 100.0 & 4h \\
 &  &  & \textbf{BB+CF+BT} & \textbf{1.39E+04} & \textbf{1} & \textbf{$\leq$0.1}$^1$ & \textbf{12} & \textbf{4.90E+03} & \textbf{107} & \textbf{$\leq$0.1}$^1$ & \textbf{16} & 1.44E+03 & 1.8E+5 & 100.0 & 4h \\
 \hline
\multirow{4}{*}{HTRU2} & \multirow{4}{*}{17,898} & \multirow{4}{*}{8} & FFT & 7.11E+04 & - & - & - & 3.36E+04 & - & - & - & 1.37E+04 & - & - & - \\
 &  &  & CPLEX & \multicolumn{4}{c}{No feasible solution} & \multicolumn{4}{c}{No feasible solution} & \multicolumn{4}{c}{No feasible solution} \\
 &  &  & BB+CF & 5.24E+04 & 1.1E+4 & $\leq$0.1 & 627 & 2.12E+04 & 2.1E+5 & 14.4 & 4h & 1.37E+04 & 9.6E+4 & 100.0 & 4h \\
 &  &  & \textbf{BB+CF+BT} & \textbf{5.24E+04} & \textbf{25} & \textbf{$\leq$0.1}$^1$ & \textbf{15} & \textbf{2.09E+04} & \textbf{5.1E+3} & \textbf{$\leq$0.1} & \textbf{282} & 1.37E+04 & 7.9E+4 & 99.5 & 4h \\
 \hline
\multirow{4}{*}{GT} & \multirow{4}{*}{36,733} & \multirow{4}{*}{11} & FFT & 4.57E+03 & - & - & - & 4.00E+03 & - & - & - & 2.59E+03 & - & - & - \\
 &  &  & CPLEX & \multicolumn{4}{c}{No feasible solution} & \multicolumn{4}{c}{No feasible solution} & \multicolumn{4}{c}{No feasible solution} \\
 &  &  & BB+CF & 3.07E+03 & 1.4E+5 & 27.6 & 4h & 2.83E+03 & 8.2E+4 & 62.2 & 4h & 2.35E+03 & 4.7E+4 & 100.0 & 4h \\
 &  &  & \textbf{BB+CF+BT} & \textbf{2.98E+03} & \textbf{2.1E+4} & \textbf{$\leq$0.1} & \textbf{2,053} & 2.81E+03 & 6.3E+4 & 60.9 & 4h & 2.29E+03 & 4.0E+4 & 100.0 & 4h \\
 \hline
\multirow{4}{*}{rds} & \multirow{4}{*}{50,000} & \multirow{4}{*}{3} & FFT & 0.11 & - & - & - & 0.06 & - & - & - & 0.03 & - & - & - \\
 &  &  & CPLEX & \multicolumn{4}{c}{No feasible solution} & \multicolumn{4}{c}{No feasible solution} & \multicolumn{4}{c}{No feasible solution} \\
 &  &  & BB+CF & 0.08 & 1.3E+5 & 4.8 & 4h & 0.05 & 8.5E+4 & 26.2 & 4h & 0.03 & 4.4E+4 & 100.0 & 4h \\
 &  &  & \textbf{BB+CF+BT} & \textbf{0.08} & \textbf{719} & \textbf{$\leq$0.1} & \textbf{16} & \textbf{0.05} & \textbf{3.6E+4} & \textbf{$\leq$0.1} & \textbf{3,429} & 0.02 & 4.1E+4 & 100.0 & 4h \\
 \hline
\multirow{4}{*}{KEGG} & \multirow{4}{*}{53,413} & \multirow{4}{*}{23} & FFT & 6.20E+06 & - & - & - & 1.70E+06 & - & - & - & 2.13E+05 & - & - & - \\
 &  &  & CPLEX & \multicolumn{4}{c}{Out of memory} & \multicolumn{4}{c}{No feasible solution} & \multicolumn{4}{c}{No feasible solution} \\
 &  &  & BB+CF & 4.98E+06 & 87 & $\leq$0.1 & 41 & 7.58E+05 & 5.8E+3 & $\leq$0.1 & 2,416 & 2.13E+05 & 2.3E+4 & 100.0 & 4h \\
 &  &  & \textbf{BB+CF+BT} & \textbf{4.98E+06} & \textbf{1} & \textbf{$\leq$0.1}$^1$ & \textbf{13} & \textbf{7.58E+05} & \textbf{1} & \textbf{$\leq$0.1}$^1$ & \textbf{14} & 2.04E+05 & 3.0E+4 & 100.0 & 4h \\
 \hline
\multirow{4}{*}{rng\_agr} & \multirow{4}{*}{199,843} & \multirow{4}{*}{7} & FFT & 4.68E+10 & - & - & - & 1.61E+10 & - & - & - & 7.47E+09 & - & - & - \\
 &  &  & CPLEX & \multicolumn{4}{c}{Out of memory} & \multicolumn{4}{c}{No feasible solution} & \multicolumn{4}{c}{No feasible solution} \\
 &  &  & BB+CF & 3.16E+10 & 3.7E+4 & 4.8 & 4h & 1.37E+10 & 2.0E+4 & 35.4 & 4h & 7.47E+09 & 1.2E+4 & 100.0 & 4h \\
 &  &  & \textbf{BB+CF+BT} & \textbf{3.14E+10} & \textbf{2.3E+3} & \textbf{$\leq$0.1}$^1$ & \textbf{239} & \textbf{1.20E+10} & \textbf{2.0E+4} & \textbf{$\leq$0.1}$^1$ & \textbf{3,330} & 7.02E+09 & 9.1E+3 & 100.0 & 4h \\
 \hline
\multirow{4}{*}{urbanGB} & \multirow{4}{*}{360,177} & \multirow{4}{*}{2} & FFT & 7.63 & - & - & - & 5.62 & - & - & - & 2.81 & - & - & - \\
 &  &  & CPLEX & \multicolumn{4}{c}{Out of memory} & \multicolumn{4}{c}{No feasible solution} & \multicolumn{4}{c}{No feasible solution} \\
 &  &  & BB+CF & 5.48 & 1.6E+4 & $\leq$0.1 & 10,713 & 4.48 & 1.5E+4 & 59.1 & 4h & 2.81 & 7.7E+3 & 100.0 & 4h \\
 &  &  & \textbf{BB+CF+BT} & \textbf{5.48} & \textbf{171} & \textbf{$\leq$0.1}$^1$ & \textbf{66} & \textbf{3.86} & \textbf{3.0E+3} & \textbf{$\leq$0.1} & \textbf{2,710} & 2.60 & 6.3E+3 & 92.6 & 4h \\
 \hline
\multirow{4}{*}{spnet3D} & \multirow{4}{*}{434,876} & \multirow{4}{*}{3} & FFT & 822.03 & - & - & - & 256.87 & - & - & - & 68.19 & - & - & - \\
 &  &  & CPLEX & \multicolumn{4}{c}{Out of memory} & \multicolumn{4}{c}{No feasible solution} & \multicolumn{4}{c}{No feasible solution} \\
 &  &  & BB+CF & 569.91 & 2.2E+4 & 0.3 & 4h & 216.25 & 1.3E+4 & 16.8 & 4h & 68.19 & 6.2E+3 & 100.0 & 4h \\
 &  &  & \textbf{BB+CF+BT} & \textbf{569.80} & \textbf{85} & \textbf{$\leq$0.1}$^1$ & \textbf{28} & \textbf{205.89} & \textbf{3.5E+3} & \textbf{$\leq$0.1}$^1$ & \textbf{661} & 68.19 & 4.6E+3 & 100.0 & 4h \\
 \hline
\end{tabular}
\begin{tablenotes}
    \item[1] Can assign $K$ initial seeds through FFT at the root node. \texttt{BB+CF+BT} results without this superscript means can not assign initial seeds.
  \end{tablenotes}
\end{threeparttable}
\end{center}
\end{table}

\begin{table}[!htp]
\vskip -0.2in
\caption{Serial results on datasets with millions of samples}
\label{tab:mil_result}
\tiny
\setlength\tabcolsep{1pt}
\begin{center}
\begin{threeparttable}
\renewcommand{\arraystretch}{1.5}
\begin{tabular}{cccccccccccccccc}
\hline
\multirow{2}{*}{Dataset} & \multirow{2}{*}{\begin{tabular}[c]{@{}c@{}}Sam\\      ple\end{tabular}} & \multirow{2}{*}{\begin{tabular}[c]{@{}c@{}}Dimen\\      sion\end{tabular}} & \multirow{2}{*}{Method} & \multicolumn{4}{c}{K=3} & \multicolumn{4}{c}{K=5} & \multicolumn{4}{c}{K=10} \\ \cline{5-16} 
 &  &  &  & UB & Nodes & \begin{tabular}[c]{@{}c@{}}Gap\\      (\%)\end{tabular} & \begin{tabular}[c]{@{}c@{}}Time\\      (s)\end{tabular} & UB & Nodes & \begin{tabular}[c]{@{}c@{}}Gap\\      (\%)\end{tabular} & \begin{tabular}[c]{@{}c@{}}Time\\      (s)\end{tabular} & UB & Nodes & \begin{tabular}[c]{@{}c@{}}Gap\\      (\%)\end{tabular} & \begin{tabular}[c]{@{}c@{}}Time\\      (s)\end{tabular} \\
 \hline
\multirow{4}{*}{USC1990} & \multirow{4}{*}{2,458,285} & \multirow{4}{*}{68} & FFT & 2.04E+11 & - & - & - & 7.47E+10 & - & - & - & 1.87E+10 & - & - & - \\
 &  &  & CPLEX & \multicolumn{4}{c}{No feasible solution} & \multicolumn{4}{c}{No feasible solution} & \multicolumn{4}{c}{No feasible solution} \\
 &  &  & BB+CF & 1.69E+11 & 916 & 3.6 & 4h & 7.47E+10 & 352 & 100.0 & 4h & 1.87E+10 & 168 & 100.0 & 4h \\
 &  &  & \textbf{BB+CF+BT} & \textbf{1.68E+11} & \textbf{1} & \textbf{$\leq$0.1}$^1$ & \textbf{277} & \textbf{6.05E+10} & \textbf{256} & \textbf{$\leq$0.1}$^1$ & \textbf{1,781} & 1.87E+10 & 396 & 61.0 & 4h \\
 \hline
\multirow{4}{*}{\begin{tabular}[c]{@{}c@{}}Gas\_\\      methane\end{tabular}} & \multirow{4}{*}{4,178,504} & \multirow{4}{*}{18} & FFT & 1.31E+08 & - & - & - & 1.17E+08 & - & - & - & 6.95E+07 & - & - & - \\
 &  &  & CPLEX & \multicolumn{4}{c}{No feasible solution} & \multicolumn{4}{c}{No feasible solution} & \multicolumn{4}{c}{No feasible solution} \\
 &  &  & BB+CF & 1.04E+08 & 1.2E+3 & 31.1 & 4h & 8.82E+07 & 488 & 100.0 & 4h & 6.95E+07 & 244 & 100.0 & 4h \\
 &  &  & \textbf{BB+CF+BT} & \textbf{1.02E+08} & \textbf{65} & \textbf{$\leq$0.1}$^1$ & \textbf{1,272} & 7.21E+07 & 807 & 19.9 & 4h & 6.95E+07 & 410 & 100.0 & 4h \\
 \hline
\multirow{4}{*}{\begin{tabular}[c]{@{}c@{}}Gas\_\\      CO\end{tabular}} & \multirow{4}{*}{4,208,261} & \multirow{4}{*}{18} & FFT & 8.83E+08 & - & - & - & 5.17E+08 & - & - & - & 2.05E+08 & - & - & - \\
 &  &  & CPLEX & \multicolumn{4}{c}{No feasible solution} & \multicolumn{4}{c}{No feasible solution} & \multicolumn{4}{c}{No feasible solution} \\
 &  &  & BB+CF & 5.66E+08 & 1.1E+3 & 12.8 & 4h & 4.09E+08 & 449 & 100.0 & 4h & 2.05E+08 & 241 & 100.0 & 4h \\
 &  &  & \textbf{BB+CF+BT} & \textbf{5.46E+08} & \textbf{66} & \textbf{$\leq$0.1} & \textbf{1,053} & \textbf{2.80E+08} & \textbf{670} & \textbf{$\leq$0.1} & \textbf{9,612} & 2.05E+08 & 398 & 100.0 & 4h \\
 \hline
\multirow{4}{*}{kddcup} & \multirow{4}{*}{4,898,431} & \multirow{4}{*}{38} & FFT & 4.71E+17 & - & - & - & 9.71E+16 & - & - & - & 5.96E+14 & - & - & - \\
 &  &  & CPLEX & \multicolumn{4}{c}{No feasible solution} & \multicolumn{4}{c}{No feasible solution} & \multicolumn{4}{c}{No feasible solution} \\
 &  &  & BB+CF & 2.25E+17 & 63 & $\leq$0.1 & 2,461 & 4.73E+16 & 229 & 100.0 & 4h & 5.96E+14 & 124 & 100.0 & 4h \\
 &  &  & \textbf{BB+CF+BT} & \textbf{2.25E+17} & \textbf{37} & \textbf{$\leq$0.1} & \textbf{958} & \textbf{4.73E+16} & \textbf{417} & \textbf{$\leq$0.1} & \textbf{10,116} & \textbf{2.58E+14} & \textbf{1} & \textbf{$\leq$0.1}$^1$ & \textbf{586} \\
 \hline
\multirow{4}{*}{HIGGS} & \multirow{4}{*}{11,000,000} & \multirow{4}{*}{29} & FFT & 368.35 & - & - & - & 320.91 & - & - & - & 198.71 & - & - & - \\
 &  &  & CPLEX & \multicolumn{4}{c}{No feasible solution} & \multicolumn{4}{c}{No feasible solution} & \multicolumn{4}{c}{No feasible solution} \\
 &  &  & BB+CF & 247.03 & 368 & 67.9 & 4h & 249.45 & 210 & 100.0 & 4h & 198.71 & 98 & 100.0 & 4h \\
 &  &  & BB+CF+BT & 237.91 & 290 & 65.5 & 4h & 235.68 & 185 & 100.0 & 4h & 198.71 & 100 & 100.0 & 4h \\
 \hline
\multirow{4}{*}{BigCross} & \multirow{4}{*}{11,620,300} & \multirow{4}{*}{56} & FFT & 1.43E+07 & - & - & - & 7.54E+06 & - & - & - & 4.24E+06 & - & - & - \\
 &  &  & CPLEX & \multicolumn{4}{c}{No feasible solution} & \multicolumn{4}{c}{No feasible solution} & \multicolumn{4}{c}{No feasible solution} \\
 &  &  & BB+CF & 1.09E+07 & 148 & 32.9 & 4h & 7.54E+06 & 122 & 100.0 & 4h & 4.24E+06 & 66 & 100.0 & 4h \\
 &  &  & BB+CF+BT & 9.97E+06 & 211 & 19.7 & 4h & 7.54E+06 & 135 & 100.0 & 4h & 4.24E+06 & 66 & 100.0 & 4h \\
 \hline
\multirow{4}{*}{\begin{tabular}[c]{@{}c@{}}Phones\\      \_acceler\\      ometer\end{tabular}} & \multirow{4}{*}{13,062,475} & \multirow{4}{*}{6} & FFT & 2.04E+28 & - & - & - & 1.09E+28 & - & - & - & 3.89E+26 & - & - & - \\
 &  &  & CPLEX & \multicolumn{4}{c}{No feasible solution} & \multicolumn{4}{c}{No feasible solution} & \multicolumn{4}{c}{No feasible solution} \\
 &  &  & BB+CF & 1.46E+28 & 51 & $\leq$0.1 & 2,038 & 6.17E+27 & 303 & 100.0 & 4h & 3.89E+26 & 148 & 100.0 & 4h \\
 &  &  & \textbf{BB+CF+BT} & \textbf{1.46E+28} & \textbf{-} & \textbf{$\leq$0.1}$^1$ & \textbf{309} & 6.17E+27 & 354 & 100.0 & 4h & 3.89E+26 & 154 & 100.0 & 4h \\
 \hline
\multirow{4}{*}{\begin{tabular}[c]{@{}c@{}}Phones\\      \_gyro\\      scope\end{tabular}} & \multirow{4}{*}{13,932,632} & \multirow{4}{*}{6} & FFT & 1.51E+28 & - & - & - & 1.09E+28 & - & - & - & 2.61E+26 & - & - & - \\
 &  &  & CPLEX & \multicolumn{4}{c}{No feasible solution} & \multicolumn{4}{c}{No feasible solution} & \multicolumn{4}{c}{No feasible solution} \\
 &  &  & BB+CF & 1.46E+28 & 51 & $\leq$0.1 & 2,195 & 6.18E+27 & 286 & 100.0 & 4h & 2.61E+26 & 136 & 100.0 & 4h \\
 &  &  & \textbf{BB+CF+BT} & \textbf{1.46E+28} & \textbf{1} & \textbf{$\leq$0.1}$^1$ & \textbf{294} & 6.18E+27 & 330 & 100.0 & 4h & 2.61E+26 & 140 & 100.0 & 4h \\
 \hline
\multirow{4}{*}{AADP} & \multirow{4}{*}{14,057,567} & \multirow{4}{*}{3} & FFT & 3.82E+03 & - & - & - & 2.98E+03 & - & - & - & 1.90E+03 & - & - & - \\
 &  &  & CPLEX & \multicolumn{4}{c}{No feasible solution} & \multicolumn{4}{c}{No feasible solution} & \multicolumn{4}{c}{No feasible solution} \\
 &  &  & BB+CF & 2.66E+03 & 602 & 35.9 & 4h & 2.49E+03 & 324 & 100.0 & 4h & 1.90E+03 & 147 & 100.0 & 4h \\
 &  &  & \textbf{BB+CF+BT} & \textbf{2.55E+03} & \textbf{196} & \textbf{$\leq$0.1} & \textbf{4,321} & 2.46E+03 & 290 & 98.9 & 4h & 1.90E+03 & 145 & 100.0 & 4h \\
 \hline
\end{tabular}
\begin{tablenotes}
    \item[1] Can assign $K$ initial seeds through FFT at the root node. \texttt{BB+CF+BT} results without this superscript means can not assign initial seeds.
\end{tablenotes}
\end{threeparttable}
\end{center}
\end{table}

\begin{table}[!htp]
\caption{Parallel results of \texttt{BB+CF+BT} ($K=3$)}
\tiny
\label{tbl: para_result}
\setlength\tabcolsep{3pt}
\begin{center}
\resizebox{0.75\textwidth}{!} {
\renewcommand{\arraystretch}{1}
\begin{tabular}{cccccccc}
\hline
Dataset & Sample & Dimension & Method & UB & Nodes & \begin{tabular}[c]{@{}c@{}}Gap\\      (\%)\end{tabular} & \begin{tabular}[c]{@{}c@{}}Time\\      (s)\end{tabular} \\
\hline
\multirow{3}{*}{HIGGS} & \multirow{3}{*}{11,000,000} & \multirow{3}{*}{29} & Heuristic & 368.35 & - & - & - \\
 &  &  & Serial & 237.91 & 290 & 65.5 & 4h \\
 &  &  & \begin{tabular}[c]{@{}c@{}}Parallel\\      (400 cores)\end{tabular} & 227.91 & 12,576 & 30.2 & 4h \\
 \hline
\multirow{3}{*}{Bigcross} & \multirow{3}{*}{11,620,300} & \multirow{3}{*}{56} & Heuristic & 1.43E+07 & - & - & - \\
 &  &  & Serial & 9.97E+06 & 211 & 19.7 & 4h \\
 &  &  & \begin{tabular}[c]{@{}c@{}}Parallel\\      (400 cores)\end{tabular} & 9.38E+06 & 10,071 & $\leq$0.1 & 6,444 \\
 \hline
\multirow{2}{*}{Taxi} & \multirow{2}{*}{1,120,841,769} & \multirow{2}{*}{12} & Heuristic & 3.09E+04 & - & - & - \\
 &  &  & \begin{tabular}[c]{@{}c@{}}Parallel\\      (2000 cores)\end{tabular} & 1.62E+04 & 1,063 & $\leq$0.1 & 5,705 \\
 \hline
\end{tabular}
}
\end{center}
\vskip -0.1in
\end{table}

\newpage

\end{document}